\documentclass[11pt,a4paper]{article}
\usepackage{amsmath}
\usepackage{amssymb}
\usepackage{amsthm}
\usepackage[all]{xy}
\usepackage{verbatim}

\begin{document}
{
\newtheorem{prop}{Proposition}[section]
\newtheorem{exa}[prop]{Example}
\newtheorem{rem}[prop]{Remark}
\newtheorem{thm}[prop]{Theorem}
\newtheorem{lem}[prop]{Lemma}
\newtheorem{de}[prop]{Definition}
\newtheorem{cor}[prop]{Corollary}

\setcounter{page}{1}
\def \Of{\mathcal{O}_{F}}
\def \Os{\mathcal{O}_{S}}
\def \Ok{\mathcal{O}_{K}}
\def \Ot{\mathcal{O}_{T}}
\def \P{\mathcal{P}}
\def \M{\mathcal{M}}
\def \N{\mathcal{N}}
\def \F{\mathcal{F}}
\def \P{\mathcal{P}}
\def \G{\mathcal{G}}
\def \pl{\varphi_{\lambda}}
\def \psil{\psi_{\lambda}}
\def \H{\mathbb{H}_{1}^{DR}}
\def \D{\mathbb{D}}
\def \t{\otimes_{\mathbb{Z}}}
\def \Osc{\mathcal{O}_{S,cris}}
\def \Md{\mathcal{M}^{DP}}
\def \Mr{\mathcal{M}^{R}}
\def \Mk{\mathcal{M}^{K}}
\def \Ms{\widetilde{\mathcal{M}}}
\def \t{\otimes_{\mathbb{Z}}}
\def \p{^{\prime}}
\def \d{^{\vee}}
\def \Qp {\mathbb{Q}_{p}}
\def \Q{\mathbb{Q}}
\def \Zp {\mathbb{Z}_{p}}
\def \Tl{T_{l}(A)}
\def \Zl{\mathbb{Z}_{l}}
\def \Zhp{\widehat{\mathbb{Z}}^{p}}
\def \Ap{\mathbb{A}^{p}_{f}}
\def \Tp{\widehat{T}^{p}(A)}
\def \Vp{\widehat{V}^{p}(A)}
\def \Gn{\widehat{\Gamma}(n)}
\def \Spec{\operatorname{Spec}}
\def \Grass{\operatorname{Grass}}
\def \Sch{\operatorname{Sch}}
\def \End{\operatorname{End}}
\def \Hom{\operatorname{Hom}}
\def \Ann{\operatorname{Ann}}
\def \Lie{\operatorname{Lie}}
\def \tr{\operatorname{tr}}
\def \Im{\operatorname{Im}}
\def \charpol{\operatorname{charpol}}
\def \Ker{\operatorname{Ker}}
\def \dim{\operatorname{dim}}
\def \ord{\operatorname{ord}}
\def \Mat{\operatorname{M}}
\def \Res{\operatorname{Res}}
\def \Def{\operatorname{Def}}
\def \Gl{\operatorname{Gl}}
\setlength{\unitlength}{0.7mm}
\def \a{ \begin{picture}(5,5)
\put(1,0){\circle*{0.5}}
\put(3,2){\circle*{0.5}}
\put(5,4){\circle*{0.5}}
\end{picture}}

\begin{center}
\begin{huge}
On the Hilbert-Blumenthal moduli problem
\footnote{to appear in: Journal of the Inst. of Math. Jussieu,
  \copyright\ Cambridge University Press}
\\[0.8cm]
\end{huge}
\begin{large}
Inken Vollaard
\footnote{
\textsc{Mathematisches Institut, Universit\"at Bonn,
  Beringstr. 1, 53115 Bonn, Germany\\
Email Address:} vollaard@math.uni-bonn.de
}\\[0.8cm]
\end{large}
\end{center}

\noindent
{\bf Abstract:}
In this article we 
compare different conditions on abelian schemes
with real multiplication which occur in the integral models of the
Hilbert-Blumenthal Shimura variety considered by Rapoport, Deligne,
Pappas  and Kottwitz. We show that the models studied by
Deligne/Pappas and Kottwitz are
isomorphic over $\Spec\mathbb{Z}_{(p)}$. We also examine the
associated local models and
prove that they are equal.

\begin{tabbing}
\textbf{Keywords:} \= Hilbert-Blumenthal varieties, bad reduction of
Shimura varieties, \\
\> abelian schemes with real multiplication, local models  
\end{tabbing}

\begin{tabbing}
\textbf{Mathematics Subject Classification (2000):} \= Primary 14G35, 11G18,
14K10 \\
\> Secondary 14K15, 11G10
\end{tabbing} 

\section*{Introduction}
Let $F$ be a totally real number field and let $\Of$ be its ring of
integers. We fix 
an arbitrary prime number $p$ and
consider abelian schemes with $\Of$-operation
$$\iota:\Of\rightarrow \End(A)$$
over $\mathbb{Z}_{(p)}$-schemes $S$ of relative dimension equal to
$g=[F:\mathbb{Q}]$. For such an abelian scheme $A$ denote by $\P(A)$ the
sheaf of $\Of$-linear, symmetric morphisms from $A$ to its dual
abelian scheme $A\d$  for the \'etale topology on $(\Sch/S)$.
Since $A$ is equipped with an $\Of$-operation, we can define the sheaf
$A\otimes_{\Of}\P(A)$ on the \'etale site, for details see Section 1.
The first goal of this paper is to compare the following three
conditions for an abelian scheme $A$ with $\Of$-operation.\\ 

\noindent
First of all, $A$ is said to be of type (DP) if it satisfies the following
equivalent conditions (Proposition 2.2). 
\begin{enumerate}
\item[(i)] The canonical morphism
$A\otimes_{\Of}\P(A)\rightarrow A\d;\ (a,\lambda)\mapsto \lambda(a)$ 
is an isomorphism.
\item[(ii)] There exists an $\Of$-linear isogeny $A\rightarrow A\d$ of degree
prime to $p$ locally for the \'etale topology on $S$.
\item[(iii)] There exists an $\Of$-linear $p$-principal polarization of $A$
locally for the \'etale topology on $S$.
\end{enumerate}
\noindent
Second, we say that $A$ is of type (R) if the Lie algebra of $A$ is
locally on $S$ a free $\Of\t\Os$-module of rank 1.\\ 

\noindent
Third, $A$ satisfies the Kottwitz determinant condition (K) if we have
$$\charpol(\iota(x);\Lie(A))=\prod_{\varphi:F\rightarrow\overline{\mathbb{Q}}}
(T-\varphi(x))$$   
for all $x\in\Of$, 
where we take the product over all embeddings of $F$ into
$\overline{\mathbb{Q}}$. \\ 

\noindent
Note that condition (R) is an open condition on the base and (DP) and
(K) are closed 
conditions. Obviously an abelian scheme of type (R) satisfies the
determinant condition. In characteristic zero all conditions are
trivial. If $p$ is unramified in $F$, all three conditions 
are equivalent ([5] Corollary 2.9).
Yu has proved in [16] Proposition 2.8 that conditions (DP) and (K) are
equivalent in the case of abelian varieties over a field and asked about the
general case ([16] Remark 2.18). We will answer this question by the
following theorem.\\

\noindent
\textbf{Theorem 1.}\emph{
The conditions (DP) and (K) are equivalent.}\\

\noindent
Let 
$\langle\cdot,\cdot\rangle$ 
be the composition of the canonical alternating $F$-bilinear form on
$F^{2}$ with the trace $tr_{F/\mathbb{Q}}$. It is an alternating
$\mathbb{Q}$-bilinear form 
such that the elements of $F$ are self-adjoint. For a $\mathbb{Q}$-linear
automorphism $g$ of $F^{2}$ denote by $g^{\ast}$ the adjoint automorphism
with respect to $\langle\cdot,\cdot\rangle$. Let $G$ be the algebraic
group over 
$\mathbb{Q}$ given by
$$G(R)=\{g\in \Gl_{F\otimes_{\mathbb{Q}}R}((F\otimes_{\mathbb{Q}}R)^{2})\mid
gg^{\ast}\in R^{\times}\}$$
for every $\mathbb{Q}$-algebra $R$. One can define different integral
models
for the Shimura variety corresponding to the group $G$ as was done by
Deligne, Pappas, Kottwitz, Rapoport and Zink. The second goal of this article 
is to compare these models.\\
 
\noindent
Let $D^{-1}\subset F$ be the inverse different of $\Of$ and
let $D^{-1}_{+}\subset F$ be the set of totally positive elements of
$D^{-1}$. Denote by $\P(A)_{S,+}$ the set of polarizations of $\P(A)_{S}$. 
For an abelian scheme $A$ with $\Of$-operation as above, a
$D^{-1}$-polarization is a homomorphism   
$\varphi:D^{-1}\rightarrow\P(A)_{S}$
of $\Of$-modules such that $\varphi(D^{-1}_{+})$ is contained in
$\P(A)_{S,+}$ and that the induced morphism
$A\otimes_{\Of}D^{-1}\rightarrow A\d$
is an isomorphism.

Fix an integer $n\geq 3$ prime to $p$. Deligne and Pappas studied in
[5] the moduli problem $\M^{DP}$ over $\Spec\mathbb{Z}_{(p)}$ defined
by the following data up to isomorphism for a
$\mathbb{Z}_{(p)}$-scheme $S$. 
\begin{enumerate}
\item An abelian scheme $A$ over $S$ with $\Of$-operation.
\item A $D^{-1}$-polarization $\varphi$ of $A$.
\item A level $n$ structure
  $\tau:A[n]\stackrel{\sim}{\longrightarrow}(\Of/n\Of)^{2}$ (Definition 3.5). 
\end{enumerate}
We will work with the definition of a level structure as in [14], which
differs from the definition in [5].
Rapoport considered in [14] the moduli problem $\M^{R}$ defined by the data
of abelian schemes of type (R), an isomorphism
$\varphi:(D^{-1},D^{-1}_{+})\stackrel{\sim}{\longrightarrow}(\P(A),\P(A)_{+})$
of 
\'etale sheaves and a level $n$ structure. 
We prove that a $D^{-1}$-polarization of an abelian scheme $A$ of type
(DP) is equivalent to an isomorphism  
$\varphi:(D^{-1},D^{-1}_{+})\stackrel{\sim}{\longrightarrow}(\P(A),\P(A)_{+})$
of \'etale sheaves (Proposition 3.3). 

Therefore, applying Theorem 1 we obtain that $\M^{R}$ is the locus of
abelian schemes of type (R) in $\M^{DP}$. 
Using the Theorem of Grothendieck-Messing it is easy to show that
$\M^{R}$ is contained in the smooth locus of $\M^{DP}$. In fact, Yu
has proved in [16] that $\M^{R}$ coincides with the smooth locus of
$\M^{DP}$. \\        

\noindent
Let $\Gn$ be the principal congruence subgroup modulo $n$
$$\Gn=\{g\in Gl_{2}(\Of\t\Zhp)\mid \det(g)\in(\Zhp)^{\times}\text{ and }
g\equiv 1 \text{ mod } n\}.$$
Another integral model for the Shimura variety corresponding to the group
$G$ was studied by Kottwitz in [8] in the unramified case and Rapoport
and Zink in the general case [15].
It is given by the following data up to ismorphism over
$\Spec\mathbb{Z}_{(p)}$. 
\begin{enumerate}
\item An abelian scheme $A\otimes_{\mathbb{Z}}\mathbb{Z}_{(p)}$ up to
  isogeny prime to $p$ with $\Of$-operation
  $\iota\otimes_{\mathbb{Z}}\mathbb{Z}_{(p)}$ which satisfies the
  determinant condition.
\item A $\mathbb{Q}$-subspace $\overline{\lambda}$ of
  $(\P(A)_{S}\t\mathbb{Z}_{(p)})\otimes_{\mathbb{Z}_{(p)}}\mathbb{Q}$ of
  dimension one which contains a $p$-principal polarization.
\item A level structure $\overline{\eta}^{p}$ of type $\Gn$ (see
  Definition 3.10). 
\end{enumerate}

\noindent
We will show that both concepts of integral models are equivalent.\\

\noindent
\textbf{Theorem 2.} \emph{The moduli problem $\M^{DP}$ is isomorphic
  to $\M^{K}$.}\\ 

\noindent
Theorem 1 can also be viewed in the light of Deligne's models of
Shimura varieties in terms of abelian varieties with weak
polarizations ([4] 4.14). More precisely,  
let $G_1$ be the group $\Res_{F/\mathbb{Q}}(\Gl_2)$.
Instead of $\M^{DP}$ or $\M^K$ one can consider the moduli problem
$\M'$ over $\Spec \mathbb{Z}_{(p)}$ given by the following data up to
isomorphism. 
\begin{enumerate}
\item An abelian scheme $A\otimes_{\mathbb{Z}}\mathbb{Z}_{(p)}$ up to
  isogeny prime to $p$ of type (DP) with $\Of$-operation
  $\iota\otimes_{\mathbb{Z}}\mathbb{Z}_{(p)}$ which satisfies the
  determinant condition.
\item A level structure $\overline{\eta}^{p}$ of type $K^p$ for a open
  compact subgroup $K^p$  of $G_1(\Ap)$ (see Definition 3.10). 
\end{enumerate}
The moduli problem $\M'$ is not representable, but by
Deligne ([4] 4.14) there is a bijection between  $\M'(\mathbb{C})$ and
the $\mathbb{C}$-valued points of the Shimura variety for $G_1$ with
respect to the compact open subgroup $K^p \Gl_2(\Of\t\mathbb{Z}_p)$ of
$G_1(\mathbb{A}_f)$.   
Theorem 1 shows that the moduli problem obtained from the definition
of $\M'$ by using only one of the conditions (DP) or (K) is equal to
$\M'$. \\ 

\noindent
To compare condition (DP) with condition (K) we will
use the local models of $\M^{DP}$ and $\M^{K}$.
Denote by $\langle\cdot,\cdot\rangle$ the canonical alternating $\Of$-bilinear
form on $\Of^{2}$.
The local model of 
 $\M^{DP}$ is defined by the functor $\N^{DP}$ over
 $\Spec\mathbb{Z}_{(p)}$ with 
\begin{align}
\N^{DP}(R)=\{\F\subset (\Of\t R)^{2}\mid \F\text{ is an }\Of\t
R\text{-submodule, locally on}&\notag\\ 
\Spec R\text{ a direct summand as an }R\text{-module and
}\F=\F^{\perp}\},&\notag 
\end{align}
where $\F^{\perp}$ is the orthogonal complement of $\F$ with respect
to $\langle\cdot,\cdot\rangle$. 
On the other hand, consider the local model $\N^{K}$ of the determinant
condition over $\Spec\mathbb{Z}_{(p)}$ defined by
\begin{align}
\N^{K}(R)=\{\F\subset (\Of\t R)^{2}\mid \F\text{ is an }\Of\t
R\text{-submodule, locally on }\Spec R\text{ a direct }&\notag\\
\text{summand as an }R\text{-module and }\F\text{ satisfies the determinant
  condition}\}.&\notag
\end{align} 
Both functors $\N^{DP}$ and $\N^{K}$ are represented by closed
 subschemes of the Grassmannian 
 $\Grass_{\mathbb{Z}_{(p)}}^{g}((\Of\otimes_{\mathbb{Z}}\mathbb{Z}_{(p)})^{2})$
 and $\N^{DP}\cap\N^{K}$ is 
 the local model of the moduli problem $\M^{K}$.\\

\noindent
\textbf{Theorem 3.}\emph{
The local models $\N^{DP}$ and $\N^{K}$ are equal over
$\Spec\mathbb{Z}_{(p)}$ as subschemes of the Grassmannian.}\\
 
\noindent
It is easy to
see that $\N^{DP}$ is contained in $\N^{K}$ (Theorem 5.6).\footnote
{Except when $p=2$. In this case, one also has to invoke the flatness
  of $\N^{DP}$([5]).} 
The generic fibres of $\N^{DP}$ and $\N^{K}$ are equal (Proposition 5.3).  
Since $\N^{K}$ is flat over $\Spec\mathbb{Z}_{p}$ ([13]), it is equal to
the flat closure of the generic fibre, hence $\N^{K}$ is contained in
$\N^{DP}$. This proves Theorem 3 abstractly. In particular, the local model
$\N^{DP}$ is flat over $\Spec\mathbb{Z}_{p}$ which was already proved
directly by 
Deligne and Pappas ([5]). Hence for $p\neq 2$ we obtain a new proof of
this last fact.
Our third goal is to give an explicit proof of Theorem 3 in case of tame
ramification.\\

\noindent
One consequence of Theorem 3 is the following result.\\

\noindent
\textbf{Corollary} \emph{
Let $A$ be an abelian scheme as above with $D^{-1}$-polarization and let
$\tilde{A}$ be a deformation of $A$ over a nilpotent thickening of the base
scheme. Then the
$D^{-1}$-polarization deforms (necessarily uniquely) if and only if $\tilde{A}$
satisfies the Kottwitz determinant condition.}\\

\noindent
Now we give a short overview of this work. 
In the first section we examine the polarization module $\P(A)$ and the
sheaf $A\otimes_{\Of}\P(A)$ and their properties.
The second section contains the study of abelian schemes of type (R), (DP)
and (K), the proof of Theorem~1 and the corollary of Theorem~3.  
The moduli problems of Deligne/Pappas and of Kottwitz type are
explained in Section 3. Section 4 contains the proof of 
Theorem 2. In the last section we give a direct proof of Theorem 3 in
the tamely ramified case.\\ 

\noindent
Finally, I want to thank everybody who helped me with this article.
First of all, I am very grateful to M. Rapoport for his introduction
into this area of mathematics
and for
his interest in my work. I want to thank G. Faltings for
pointing out an inaccuracy in an earlier version of this work and the
referee and C.-F. Yu 
for useful comments on this article. 
Furthermore, many thanks are due to 
T. Wedhorn and U. G\"ortz for their
support and lots of helpful discussions.


\section{Notations and first properties}

Let $F$ be a totally real extension of $\mathbb{Q}$ of degree $g$. Denote
by $\Of$ the ring of integers of $F$ and fix a
prime number $p$. In this article we will only consider abelian
schemes of relative dimension $g$ 
over a locally noetherian $\Spec\mathbb{Z}_{(p)}$-scheme.
\begin{de}
A homomorphism of rings
$$\iota:\Of\rightarrow\End(A)$$
is called an $\Of$-operation of $A$. 
\end{de}

\noindent
\textbf{1.2.} 
Let $A$ be an abelian scheme with $\Of$-operation and let $A\d$ be the dual
abelian scheme. We obtain an $\Of$-operation on $A\d$ via
duality. Denote by $\P(A)$ the sheaf for 
the \'etale topology (big \'etale site) on $(\Sch/S)$ defined by
$$\P(A)_{T}=\{\lambda :A_{T}\rightarrow A_{T}\d\mid \lambda\text{ is
  symmetric and }\iota_{T}\d\circ\lambda=\lambda\circ\iota_{T}\}$$
for all $T\rightarrow S$. Let $\P(A)_{T,+}$ be the subset of $\P(A)_T$
  of all $\Of$-linear polarizations.

Every element $\lambda\in\P(A)$ not equal to zero is an
isogeny. Indeed, it is sufficient to prove this in the case of abelian
varieties. But every abelian variety of dimension $g=[F:\mathbb{Q}]$ with
$\Of$-operation is simple in the category of abelian varieties with
$\Of$-operation (compare [10]). Thus $\lambda\neq 0$ is an isogeny. 

\addtocounter{prop}{1}
\begin{prop}
Let $A$ be an abelian scheme with $\Of$-operation over a connected base
scheme $S$ and let $s$ be a geometric point of $S$.
\begin{enumerate} 
\item[a)] The $\Of$-module $\P(A)_{S}$ is trivial or an invertible
  $\Of$-module such that $\P(A)_{S,+}$
is a positivity of
$\P(A)_{S}$.
\item[b)] If $\P(A)_{S}$ is non-trivial, then
$\P(A)_{S}\otimes_{\mathbb{Z}}\mathbb{Z}_{l}=
\P(A)_{s}\otimes_{\mathbb{Z}}\mathbb{Z}_{l}$  
for every prime $l$ which is invertible on $S$.
\item[c)] Let $T\rightarrow S$ be a faithfully flat morphism and let
  $T$ be a connected scheme, 
then $\P(A)_{T}=\P(A)_{S}$ if $\P(A)_{S}$ is non-trivial. In particular, if
$S$ is the spectrum of a field $k$, the sheaf for the \'etale topology
$\P(A)$ is constant, 
i.e., the constant sheaf for the big \'etale site, associated to the
module $\P(A)_{k}$. Moreover, $\P(A)_k$ is nontrivial. 
\end{enumerate}
\end{prop} 

\noindent
Before proving Proposition 1.3 we will first use it to prove the
following result. 

\begin{prop}
Let $A$ be an abelian scheme with $\Of$-operation over a scheme $S$. If
$S$ is a $\mathbb{Q}$-scheme, the sheaf for the \'etale
topology $\P(A)$ is a locally constant sheaf of invertible
$\Of$-modules. The same is true if there exists an $\Of$-linear 
$p$-principal polarization locally for the \'etale topology.
\end{prop}

\noindent
\textsc{Proof of Proposition 1.4:}
We will first prove the proposition in the case of characteristic
zero. Therefore, we have to prove the following claim.

\emph{
Claim:} In characteristic zero the module $\P(A)$ is locally for the \'etale
topology on $S$ non
trivial.
  
If $A=\Lie(A)/H_{1}(A,\mathbb{Z})$ is an abelian variety over $\mathbb{C}$, an
$\Of$-linear polarization can be constructed by the following arguments
given by Goren in
[7] Section 2.2.
Let $\sigma_{1},...,\sigma_{g}$ be the embeddings of $F$
into $\mathbb{R}$, then 
$F\otimes_{\mathbb{Q}}\mathbb{C}$ is equal to
$\prod_{i=1}^{g}\mathbb{C}$. Since $\Lie(A)$ is an
$F\otimes_{\Q}\mathbb{C}$-module, we obtain a decomposition
$$\Lie(A)=\bigoplus_{i=1}^{g}L_{i},$$
where $L_{i}$ is a $\mathbb{C}$-vector space of dimension 1 and the
$\Of$-multiplication is given by
$$x(v_{1},...,v_{g})=(\sigma_{1}(x)v_{1},...,\sigma_{g}(x)v_{g})$$ 
for an element $x\in\Of$. 
The
lattice $H_{1}(A,\mathbb{Z})$ is $\Of$-invariant, hence equal to
$\mathfrak{a}e\oplus\mathfrak{b}f$ for two ideals $\mathfrak{a},\mathfrak{b}$
of $\Of$ and $e, f\in \Lie(A)$.  
We can define a Riemann form on $\Lie(A)$ by 
$$\psi(v,w)=\sum_{i=1}^{g}\frac{v_{i}\overline{w}_{i}}{\Im(e_{i}
  \overline{f}_{i})}.$$ 
The elements of $\Of$ are
self-adjoint with respect to $\psi$, hence $\psi$ defines an $\Of$-linear
polarization on $A$. The same proof works if $A$ is an abelian scheme over
a simply connected analytic space $U$ over $\mathbb{C}$. 

If $A$ is an abelian scheme over an arbitrary $\mathbb{C}$-scheme
$S$, we may assume that $S$ is of finite type over $\mathbb{C}$.
It is sufficient to
construct an $\Of$-linear polarization locally for the \'etale topology at
every closed point $s$ of $S$. Since there exists a simply connected analytic
neighbourhood of $s$, we obtain an $\Of$-linear polarization over the
analytic local 
ring $\mathcal{O}_{S,s}^{an}$, hence over the completion
$\widehat{\mathcal{O}}_{S,s}$ of the algebraic local ring. Using Artin's
approximation theorem ([1] Theorem 1.12) we obtain an $\Of$-linear
polarization over the 
henselization of $\mathcal{O}_{S,s}$ and the claim is proved in this case.

Now let $S$ be an arbitrary scheme of finite type over
$\mathbb{Q}$. Using the above result  
we obtain an $\Of$-linear polarization locally for the
\'etale topology over the base change $S_{\mathbb{C}}$. 
Since the claim is local on $S$,
we may assume
  that $S$ is affine and that there exists a surjective \'etale morphism
  $U\rightarrow S_{\mathbb{C}}$ together with a polarization over
  $U$. The morphism 
  $U\rightarrow S$ is faithfully flat and has regular fibres. Using N\'eron
  desingularization ([12] Theorem 2.5) we obtain a scheme $V$ with a
  smooth and surjective morphism 
   $f:V\rightarrow S$ such that there exists an $\Of$-linear
   polarization over $V$.  
Since $f$ is smooth and surjective, there exists a quasi-section of $f$
which means that there exists an \'etale scheme $S\p\twoheadrightarrow S$
and an $S$-morphism $S\p\rightarrow V$ ([6] 17.16.3). We obtain an
$\Of$-linear polarization on $A_{S\p}$ by base change from $V$ to
$S\p$  which proves the claim. 

To finish the proof of the proposition in characteristic zero we may
assume that $S$ is connected 
and $\P(A)_{S}$ is non-trivial.
We obtain from Proposition 1.3 that $\P(A)_{S}$ is isomorphic to
$\P(A)_{s}$. The same is true for every connected scheme $T$ over
$S$. Let $t$ be a geometric point of $T$ and let $s$ be the underlying
geometric 
point of $S$. Since $t\rightarrow s$ is faithfully flat, the claim is proved
using part c) of Proposition 1.3. 

If there exists an $\Of$-linear $p$-principal polarization, use the
fact that an $\Of$-linear $p$-principal polarization over $S$
generates 
the $\Of\otimes_{\mathbb{Z}}\mathbb{Z}_{p}$-module
$\P(A)_{S}\otimes_{\mathbb{Z}}\mathbb{Z}_{p}$.\hfill{$\Box$}\\

\noindent
\textbf{1.5.} 
Denote by $D^{-1}$ the inverse ideal of the different of $\Of$. Let $R$ be
a ring and let $M$ be an $\Of\otimes_{\mathbb{Z}}R$-module. We obtain
an isomorphism 
\begin{align}
\Hom_{R\otimes_{\mathbb{Z}}\Of}(M,R\otimes_{\mathbb{Z}}D^{-1})&
\stackrel{\sim}{\longrightarrow}\Hom_{R}(M,R)\tag{1.5.1}\\  
f&\longmapsto \tr_{F/\mathbb{Q}}\circ f.\notag
\end{align}
Thus choosing an $R$-bilinear form $\Phi':M\times M\rightarrow R$ such
that the elements of $\Of$ are self-adjoint, is the same as choosing
an $R\t\Of$-bilinear form 
$\Phi:M\times M\rightarrow R\t D^{-1}$.
\addtocounter{prop}{1}
\begin{prop}
Let $A$ be an abelian variety with
$\Of$-operation over an algebraically closed field $k$. Then the Tate
module $\Tl$ is a free $\Of\t\Zl$-module of 
rank 2 for every prime number $l\neq \text{char } k$.
\end{prop}

\noindent 
\textsc{Proof:} The proof is analogous to the proof of [14] Lemma
1.3.\hfill{$\Box$}\\ 

\noindent
\textbf{1.7.}
\textsc{Proof of Proposition 1.3.:}\\ 
Let $\pi_{1}(S,s)$ be the fundamental group.
The morphism
\begin{align}
\Hom(A,A\d)\hookrightarrow\Hom_{\mathbb{Z}_{l}}(T_{l}(A),
T_{l}(A\d))\tag{1.7.1}
\end{align}
is injective because the union of the schemes $A[l^{n}]$ is scheme
theoretically dense 
in $A$. Hence (1.7.1) induces an injective morphism
$$\P(A)_{S}\otimes_{\mathbb{Z}}\mathbb{Z}_{l}\hookrightarrow
\Hom_{\Of\otimes_{\mathbb{Z}}\mathbb{Z}_{l}}(T_{l}(A), 
T_{l}(A\d)).$$ 
Denote by $\mu_{l^{\infty}}=\underleftarrow{\lim}\ \mu_{l^{r}}$
the inverse limit of the
$l^r$th roots of unity for all positive integers $r$ and let
$\lambda$ be an element of  
$\P(A)_{S}$. 
Using the canonical pairing 
$$e_{l}:T_{l}(A)\times T_{l}(A\d)\rightarrow\mu_{l^{\infty}}$$
the morphism $\lambda$ induces an alternating form 
$\psi'_\lambda$ on $T_{l}(A)$. Since $\lambda$ is $\Of$-linear, the
Rosati involution with respect to $\lambda$ is trivial on $\Of$. Hence
the elements of $\Of$ are self-adjoint with respect to
$\psi'_\lambda$. 
By (1.5.1) we obtain from $\psi'_\lambda$ an alternating
$\Of\otimes\mathbb{Z}_{l}$-bilinear form
$$\psi_{\lambda}:T_{l}(A)\times T_{l}(A)\rightarrow
D^{-1}\otimes_{\mathbb{Z}}\mu_{l^{\infty}}$$ 
such that $tr_{F/\mathbb{Q}}\circ \psi_{\lambda}=e_{l}(\cdot,\lambda(\cdot))$.
Thus we obtain injective morphisms
\begin{align}
\P(A)_{S}\otimes_{\mathbb{Z}}\mathbb{Z}_{l}&\hookrightarrow
\Hom_{\Of\otimes_{\mathbb{Z}}\mathbb{Z}_{l},alt}(T_{l}(A_{s}),  
T_{l}(A_{s}\d))^{\pi_{1}(S,s)}\tag{1.7.2}\\
&\hookrightarrow\Hom_{\Of\otimes_{\mathbb{Z}}\mathbb{Z}_{l},alt}(T_{l}(A_{s}),
T_{l}(A_{s}\d))\tag{1.7.3}\\
&=\Hom_{\Of\otimes_{\mathbb{Z}}\mathbb{Z}_{l}}(\wedge_{\Of\otimes_{\mathbb{Z}}
  \mathbb{Z}_{l}}^{2}T_{l}(A_{s}),D^{-1}\otimes_{\mathbb{Z}}\mu_{l^{\infty}}). 
\notag    
\end{align}
The $\Of$-module $\P(A)_{S}$ is locally free because it is torsion
free and of finite type, 
hence $\P(A)_{S}\otimes_{\mathbb{Z}}\mathbb{Z}_{l}$ is a free
$\Of\otimes_{\mathbb{Z}}\mathbb{Z}_{l}$-module. On 
the other hand, the Tate module $T_{l}(A_{s})$ is a free
$\Of\otimes_{\mathbb{Z}}\mathbb{Z}_{l}$-module  of rank 2 (Proposition
1.6). Therefore, we obtain an isomorphism between 
$\Hom_{\Of\otimes_{\mathbb{Z}}\mathbb{Z}_{l}}(\wedge_{\Of\otimes_{\mathbb{Z}}
  \mathbb{Z}_{l}}^{2}T_{l}(A_{s}), 
D^{-1}\otimes_{\mathbb{Z}}\mu_{l^{\infty}})$ and
$\Of\otimes_{\mathbb{Z}}\mathbb{Z}_{l}$. Since both morphisms (1.7.2)
and (1.7.3) have torsion free cokernels, the module $\P(A)_{S}$ is
equal to zero or an invertible $\Of$-module. In the second case the
morphism  
\begin{align}
\P(A)_{S}\otimes_{\mathbb{Z}}\mathbb{Z}_{l}
\stackrel{\sim}{\longrightarrow}\Hom_{\Of\otimes_{\mathbb{Z}}
  \mathbb{Z}_{l},alt}(T_{l}(A_{s}), 
T_{l}(A_{s}\d))\tag{1.7.4}
\end{align}
is an isomorphism. The set $\P(A)_{S,+}$ is a positivity of $\P(A)_{S}$
([14] Proposition 1.17), hence a) is proved. 

To prove b) note that for every abelian variety there exists an
$\Of$-linear polarization ([14] Proposition 1.12). Hence $\P(A)_{s}$ is non
trivial and the same proof as above shows that
$\P(A)_{s}\otimes_{\mathbb{Z}}\mathbb{Z}_{l}$
is equal to $\Hom_{\Of\otimes_{\mathbb{Z}}\mathbb{Z}_{l},alt}(T_{l}(A_{s}),
T_{l}(A_{s}\d))$.

To prove c) note that the morphism of invertible $\Of$-modules
$f:\P(A)_{S}\rightarrow\P(A)_{T}$ is 
injective and that the cokernel is $\mathbb{Z}$-torsion. To prove the
surjectivity let $\lambda$ be an element of 
$\P(A)_{T}$. There exists an integer $n$ such that $n\lambda=\mu_{T}$
for a morphism $\mu\in\P(A)_{S}$.
The morphism $n\lambda$ vanishes on $A_{T}[n]$. 
As $T\rightarrow S$ is faithfully flat, the morphism $\mu$ vanishes on
$A[n]$. Thus 
there exists a morphism $\lambda\p\in\P(A)_{S}$ such that
$\lambda\p_{T}=\lambda$, hence $f$ is an isomorphism.\hfill{$\Box$}\\

\noindent
\textbf{1.8.} Let $A$ be an abelian scheme with $\Of$-operation over $S$
and let $L$ be an invertible $\Of$-module.
Denote by $A\otimes_{\Of}L$ the sheaf for the \'etale topology on
$(\Sch/S)$ associated to the presheaf of $\Of$-modules
$$T\mapsto A(T)\otimes_{\Of} L.$$ 
Consider an element $\lambda$ of $L$ and denote by $\mathfrak{a}_{\lambda}$ the
annihilator ideal $\Ann_{\Of}(L/\lambda\Of)$ of $\Of$.
Then the sequence 
\begin{align}
0\rightarrow A[\mathfrak{a}_{\lambda}]\rightarrow A&\rightarrow
A\otimes_{\Of}L \rightarrow 0\tag{1.8.1}\\
a&\mapsto a\otimes\lambda\notag
\end{align}
is an exact sequence of \'etale sheaves of $\Of$-modules over $S$. The
rank of $A[\mathfrak{a}_{\lambda}]$ is equal to
$N(\mathfrak{a})^{2}=\mid\!L/\lambda\Of\!\mid^{2}$. We can 
choose an element $\lambda$ such that its degree is prime to
$p$. Then $A[\mathfrak{a}_{\lambda}]$ is a finite, \'etale group scheme
and the quotient $A\otimes_{\Of}L$ for the \'etale topology is represented by a
scheme. 
Let $A\otimes_{\Of}\P(A)$ be the tensor product of $A$ with the sheaf for
the \'etale topology $\P(A)$ and 
denote by $\tau$ the morphism of sheaves for the \'etale topology 
\begin{align}
\tau:A\otimes_{\Of}\P(A)&\rightarrow A\d\notag\\
(a,\lambda)&\mapsto \lambda(a).\notag
\end{align}
The sheaf $A\otimes_{\Of}\P(A)$ is not representable by a scheme in
general (Example 2.16) but obviously it is representable if $\P(A)$
is locally constant. The converse is true as well.

\addtocounter{prop}{2} 
\begin{lem}
The sheaf for the \'etale topology $A\otimes_{\Of}\P(A)$ is
represented by an abelian scheme 
if and only if $\P(A)$ is locally constant.
\end{lem}

\noindent
\textsc{Proof:}
Let $A\otimes_{\Of}\P(A)$ be represented by a scheme. The
module $\P(A)$ is locally on $S$ non-trivial. We may assume that $S$
is the spectrum 
of a local henselian ring. Let $P$ be the invertible module
$\P(A)_{S}$ and let $K$ be 
the kernel of the isogeny 
$$\mu:A\otimes_{\Of}P\rightarrow A\otimes_{\Of}\P(A).$$
Since $\P(A)_{T}$ is equal to $\P(A)_{S}$ for every connected scheme $T$
faithfully flat
over $S$ (Proposition 1.3), the kernel $K(T)$ is trivial for all $T$
faithfully flat over $S$. Therefore, $K(K)$ is equal to zero which means that
$K$ is trivial. We obtain that $\mu$ is an isomorphism. 

Let $s$ be a
geometric point of $S$. The sheaf $\P(A_{s})$ is the constant sheaf
associated to the module $P\p=\P(A)_{s}$ (Proposition 1.3).
For $\P(A)$ to be constant it is sufficient to show that the morphism 
$$\varphi:P\rightarrow P\p$$
is an isomorphism. 
It is clear that $\varphi$
is injective. We can prove the surjectivity of $\varphi$ after localizing
at any prime ideal of $\Of$. Hence we can assume that $P$
and $P\p$ are free $\Of$-modules. Let $\lambda$ be a generator of $P\p$
and let $x\lambda$ be a generator of $P$,
where $x$ is an element of $\Of$. 
The isomorphism $\mu$ induces an isomorphism
$$A_{s}\otimes_{\Of}(x\lambda)\stackrel{\sim}{\longrightarrow}
A_{s}\otimes_{\Of}(\lambda),$$
hence
its kernel $A_{s}[x]$ is equal to zero. But the degree of the
multiplication with $x$ on $A_{s}$ is equal to $N(x)^{2}$. Therefore,
$x$ is an element of $\Of^{\times}$ and $\varphi$ is an 
isomorphism.\hfill{$\Box$}\\

\begin{prop}
Let $A$ be an abelian scheme over a connected scheme $S$ with
$\Of$-operation and let $l$ be a prime
number which is invertible on $S$. Suppose that $P:=\P(A)_{S}$ is non
trivial. Then the isogeny  
\begin{align}
\tilde{\tau}:A\otimes_{\Of}P&\rightarrow A\d\notag\\
(a,\lambda)&\mapsto \lambda(a)\notag
\end{align}
is of degree prime to $l$. In particular, if the polarization module
$\P(A)$ is locally constant, 
the morphism
$\tau:A\otimes_{\Of}\P(A)\rightarrow A\d$ is an isogeny of degree prime to $l$.
\end{prop} 

\noindent
\textsc{Proof:} We will show that for every geometric point $s$ of $S$
the morphism $\tilde{\tau}_{s}$ induces an isomorphism on the $l$-adic
Tate modules. Then the kernel of $\tilde{\tau}$ will be prime to $l$. 
The Tate module of 
$A_{s}\otimes_{\Of}P$ is equal to $T_{l}(A_{s})\otimes_{\Of}P$, hence
$\tilde{\tau}_{s}$ induces a morphism 
$$T_{l}(\tilde{\tau}_{s}):T_{l}(A_{s})\otimes_{\Of}P\rightarrow
T_{l}(A_{s}\d).$$ 
Since
$P\otimes_{\mathbb{Z}}\mathbb{Z}_{l}\stackrel{\sim}{\longrightarrow}
\Hom_{\Of\otimes_{\mathbb{Z}}\mathbb{Z}_{l},alt}(T_{l}(A_{s}),T_{l}(A_{s}\d))$ 
is an isomorphism (1.7.4), we obtain an isomorphism
$$ (\wedge_{\Of\t\Zl}^{2}T_{l}(A_{s}))\cong D^{-1}P^{-1}\t\Zl.$$
This induces a perfect $\Zl$-bilinear form
$$T_{l}(A_{s})\times T_{l}(A_{s})\rightarrow P^{-1}\otimes_{\Of}\Zl$$
which is equal to the $\Zl$-bilinear form induced by $T_{l}(\tilde{\tau}_{s})$.
Hence $T_{l}(\tilde{\tau}_{s})$ is an isomorphism.

\hfill{$\Box$}

\noindent
\begin{cor}
Let $A$ be an abelian scheme with $\Of$-operation over a
$\mathbb{Q}$-scheme $S$. The polarization module $\P(A)$ is locally
constant and the canonical isogeny
$\tau:A\otimes_{\mathbb{Z}}\P(A)\rightarrow A\d$ 
is an isomorphism.
\end{cor}

\noindent
\textsc{Proof:} The corollary follows from Proposition 1.4 and
Proposition 1.10.\hfill{$\Box$}\\ 

\noindent
\textbf{Remark 1.12.}
Let $A$ be an abelian scheme with $\Of$-operation over a
$\mathbb{Q}$-scheme $S$. Then the Lie algebra of $A$ is locally on $S$ a
free $\Of\otimes_{\mathbb{Z}}\mathcal{O}_{S}$-module. Indeed,
for every
field $k$ of characteristic zero the operation of
$\Of\otimes_{\mathbb{Z}}k$ on the Lie algebra is faithful. Hence the Lie
algebra is a free $\Of\otimes_{\mathbb{Z}}k$-module ([14] Proposition 1.4).\\

\noindent
\textbf{Remark 1.13.}
There exist abelian varieties with $\Of$-operation in characteristic
$p$ which cannot be lifted 
to characteristic zero, see [14] 1.29 for an example in the
case when $p$ is unramified in $\Of$. Therefore, to define a flat
moduli problem we have to add some 
conditions in characteristic $p$.


\section{Abelian schemes with real multiplication}
In this section we will compare different properties of abelian
schemes with real multiplication. 

\begin{de} 
Let $A$ be an abelian scheme with $\Of$-operation over $S$. Then $A$
is of type (R) 
if the Lie algebra of $A$ is locally on $S$ a free
$\Of\otimes_{\mathbb{Z}}\mathcal{O}_{S}$-module (see [14]).
\end{de}

\begin{prop}
Let $A$ be an abelian scheme over $S$ with $\Of$-operation. Then the
following properties are equivalent. 
\begin{itemize}
\item[a)] The canonical morphism $\tau:A\otimes_{\Of}\P(A)\rightarrow A\d;\
  (a,\lambda)\mapsto\lambda(a)$ is an isomorphism.
\item[b)] There exists an $\Of$-linear $p$-principal polarization on $A$
  locally for the \'etale topology on $S$.
\item[c)] There exists an isogeny $A\rightarrow A\d$ of degree prime
  to $p$ locally for the 
  \'etale topology on $S$.
\end{itemize}
The $\Of$-linear $p$-principal polarization is uniquely determined up
to a totally 
positive element $x$ of $(\Of\t\mathbb{Z}_{(p)})^{\times}$ such that
$x\lambda$ is a morphism.
\end{prop}

\begin{de} 
We call an abelian scheme $A$ with $\Of$-operation of type (DP) if it
satisfies the equivalent conditions of Proposition 2.2 (see [5]).
\end{de}

\begin{de}
An abelian scheme $A$ with $\Of$-operation satisfies the Kottwitz
determinant condition (K) if we have
$$\charpol(x;\Lie(A))=\prod_{\varphi:F\rightarrow
  \mathbb{\overline{Q}}}(T-\varphi(x))$$ 
for  every element $x\in\Of$ (see [8] Section 5).
\end{de}

\noindent
\textbf{Remark 2.5.}
Condition (R) is an open condition and the conditions (DP) and (K) are
closed conditions.
Every abelian scheme over
a $\mathbb{Q}$-base scheme satisfies conditions (R), (K) and (DP)
(Corollary 1.11, Remark 1.12),
thus all conditions are trivial in characteristic zero. Obviously
condition (R) implies the determinant condition (K). We will show 
that conditions (DP) and (K) are equivalent. This was already proved by Yu
([16] Proposition 2.8)
in the case of abelian varieties.
Furthermore, all conditions are equivalent if $p$ is
unramified ([5] Corollary 2.9; cf. Remark 5.5). \\

\noindent
\textbf{2.6.} Let $A$ be an abelian variety over an algebraically closed
field $k$ of characteristic $p$. Then the determinant condition can be
expressed in the following way as proved by Yu in [16]. The decomposition
$$\Of\otimes_{\mathbb{Z}}k=\prod_{\mathfrak{p}|p}
\widehat{\mathcal{O}}_{F,\mathfrak{p}}\otimes_{\mathbb{Z}_{p}}k$$   
induces a decomposition of the Lie algebra
$$\Lie(A)=\bigoplus_{\mathfrak{p}|p}\Lie(A)_{\mathfrak{p}}.$$ 
Let $W(k)$ be the ring of Witt vectors and
denote by $\widehat{\mathcal{O}}_{F,\mathfrak{p}}^{ur}$ the maximal unramified
extension of $\mathbb{Z}_{p}$ in $\widehat{\mathcal{O}}_{F,\mathfrak{p}}$. We
obtain 
$$\widehat{\mathcal{O}}_{F,\mathfrak{p}}\otimes_{\mathbb{Z}_{p}}k=
\prod_{i=1}^{f_{\mathfrak{p}}}k[\pi]/(\pi^{e_{\mathfrak{p}}})$$  
for a
uniformizing element $\pi$ in $\widehat{\mathcal{O}}_{F,\mathfrak{p}}$.
We obtain a decomposition
$$\Lie(A)_{\mathfrak{p}}=\bigoplus_{i=1}^{f_{\mathfrak{p}}}
L_{\mathfrak{p},i},$$  
where $L_{\mathfrak{p},i}$ is a
$k[\pi]/(\pi^{e_{\mathfrak{p}}})$-module and the elements of 
$\widehat{\mathcal{O}}_{F,\mathfrak{p}}^{ur}$ operate on
$L_{\mathfrak{p},i}$ via the different 
embeddings 
$\sigma_{\mathfrak{p},i}:\widehat{\mathcal{O}}_{F,\mathfrak{p}}\rightarrow
W(k)$. The Lie algebra $\Lie(A)$ is a quotient of the free $\Of\t
k$-module $\H(A/k)$ of rank 2. Using the analogous decomposition of
the de~Rham homology $\H(A/k)$ of $A$ 
$$\H(A/k)=\bigoplus_{\mathfrak{p}|p}\bigoplus_{i=1}^{f_{\mathfrak{p}}}
H_{\mathfrak{p},i}$$ 
into free $k[[\pi]]/(\pi^{e_{\mathfrak{p}}})$-modules
$H_{\mathfrak{p},i}$ of rank 2, the
$k[[\pi]]/(\pi^{e_{\mathfrak{p}}})$-module $L_{\mathfrak{p},i}$ is a
quotient of $H_{\mathfrak{p},i}$ for every $i$ and
$\mathfrak{p}$. Hence we obtain a decomposition  
$$L_{\mathfrak{p},i}=k[\pi]/(\pi^{e_{1}^{\mathfrak{p},i}})\oplus
k[\pi]/(\pi^{e_{2}^{\mathfrak{p},i}})$$ 
with integers $0\leq e_{1}^{\mathfrak{p},i},e_{2}^{\mathfrak{p},i}\leq
e_{\mathfrak{p}}$ such that 
$\sum_{\mathfrak{p}\mid
  p}\sum_{i=1}^{f_{\mathfrak{p}}}(e_{1}^{\mathfrak{p},i}+
e_{2}^{\mathfrak{p},i})=g$.   

Then $\Lie(A)$ satisfies the determinant condition if and only if 
$e_{1}^{\mathfrak{p},i}+e_{2}^{\mathfrak{p},i}=e_{\mathfrak{p}}$ for all
$\mathfrak{p}$ and $i$.
In particular,
if $p$ is totally ramified in $\Of$, the
determinant condition is trivial. 
If $p$ is unramified, we obtain a decomposition
$\Lie(A)=\bigoplus_{i=1}^{g}L_{i}$, 
where $L_{i}$ is a $k$-vector space of dimension $\dim_{k}L_{i}\leq
2$. Then the 
determinant condition is equivalent to $\dim_{k}L_{i}=1$ for all $i$,
i.e., $\Lie(A)$ 
is a free $\Of\otimes_{\mathbb{Z}}k$-module of rank 1.\\

\noindent
\textbf{2.7.}
\textsc{Proof of Proposition 2.2:}
The polarization module $\P(A)$ is locally constant under any of the hypotheses
(Proposition 1.4, Lemma 1.9). 
Since all conditions are local for the \'etale topology, we may assume that
$\P(A)$ is constant and $S$ is connected.
Let $\lambda\in\P(A)_{S}$ be an isogeny and consider the exact diagram of
sheaves for the \'etale topology
$$\xymatrix{
0\ar@{->}[r]  & A[\mathfrak{a}_{\lambda}]\ar@{->}[r]\ar@{=}[d]
&\Ker\lambda\ar@{->}[r]\ar@{->}[d] &\Ker\tau\ar@{->}[r]\ar@{->}[d] &0\\ 
0 \ar@{->}[r] & A[\mathfrak{a}_{\lambda}]\ar@{->}[r] & A
\ar@{->}[r]\ar@{->}[d]^\lambda
&A\otimes_{\Of}\P(A)_{S}\ar@{->}[r]\ar@{->}[d]^\tau & 0\\   
 & &A\d \ar@{=}[r] &A\d,  &\\
}$$
where $\mathfrak{a}_{\lambda}$ is the annihilator ideal of
$\P(A)_{S}/\lambda\Of$.

Let $\tau$ be an isomorphism. Choose an element 
$\lambda\neq 0$ of $\P(A)_{S}$ such that
$\mid\!\P(A)_{S}/\lambda\Of\!\mid$ is prime to $p$. Using the above diagram we
obtain that $\lambda$ is an isogeny with kernel equal to
$A[\mathfrak{a}_{\lambda}]$. But the degree of $\lambda$ is equal to
$\mid\!A[\mathfrak{a}_{\lambda}]\!\mid=\mid\!(\P(A)_{S}/\lambda\Of)\!\mid^{2}$,
hence 
prime to $p$. Now choose an element $x\in\Of$ with $N(x)$ prime to $p$
such that 
$x\lambda$ is totally positive in $\P(A)_{S}$. Then $x\lambda$ is a
$p$-principal 
polarization. Thus b) follows from a).
It is clear that b) implies c).

Now let $\lambda$ be an isogeny prime to $p$, then the morphism $\tau$
is an isogeny prime to 
$p$. On the other hand, the degree of $\tau$ is prime to $l$ for every
prime $l\neq p$ (Lemma 1.9), hence $\tau$ is an isomorphism and a)
follows from c).

The uniqueness of the $\Of$-linear $p$-principal polarization follows
from the fact that 
$x\lambda$ is a polarization for an element
$x\in(\Of\t\mathbb{Z}_{(p)})^{\times}$ if and only if $x$ is totally
positive ([17] 
Lemma 1.4). \hfill{$\Box$}\\

\noindent
\textbf{2.8.}
To show that the determinant condition is equivalent to condition (DP) we
need the following result on local models. Denote by
$\langle\cdot,\cdot\rangle$ 
the canonical alternating $\Of$-linear form on $\Of^{2}$. The local
model $\N^{DP}$ is defined by
\begin{align}
\N^{DP}(R)=\{\F\subset (\Of\t R)^{2}\mid \F\text{ is an }\Of\t
R\text{-submodule, locally on}&\notag\\ 
\Spec R\text{ a direct summand as an }R\text{-module and
}\F=\F^{\perp}\}& \notag
\end{align}
for every $\mathbb{Z}_{(p)}$-algebra $R$,
where $\F^{\perp}$ is the orthogonal complement of $\F$ with respect to
$\langle\cdot,\cdot\rangle$. Furthermore, denote by $\N^{K}$ the functor
\begin{align}
\N^{K}(R)=\{\F\subset (\Of\t R)^{2}\mid \F\text{ is an }\Of\t
R\text{-submodule, locally on }\Spec R\text{ a direct }&\notag\\
\text{summand as an }R\text{-module and }\F\text{ satisfies the determinant
  condition}\}&\notag
\end{align} 
for every $\mathbb{Z}_{(p)}$-algebra $R$. Both $\N^{DP}$ and $\N^{K}$ are
closed subschemes of the Grassmannian
$\Grass_{\mathbb{Z}_{(p)}}^{g}((\Of\otimes_{\mathbb{Z}}
\mathbb{Z}_{(p)})^{2})$.  
In Section 5 we will examine $\N^{DP}$ and $\N^{K}$ and prove the following
theorem.

\addtocounter{prop}{4}
\begin{thm}
The local models $\N^{DP}$ and $\N^{K}$ are equal over
$\Spec\mathbb{Z}_{(p)}$ as subschemes of the Grassmannian.
\end{thm}

\noindent
\textsc{Proof:} Theorem 5.6.\hfill{$\Box$}

\begin{cor}
Let $A$ be an abelian scheme over a $\mathbb{Z}_{(p)}$-scheme $S$ with
$\Of$-operation $\iota$ and $\Of$-linear polarization $\lambda$. 
Denote by $S\hookrightarrow S\p$ a thickening with locally nilpotent
divided powers and by $A\p$ a
lifting of $(A,\iota)$ to $S\p$. If $\Lie(A\p)$ satisfies
the Kottwitz determinant condition, the polarization $\lambda$ can be
lifted to an $\Of$-linear polarization of $A^{\prime}$. 
\end{cor}

\noindent
\textsc{Proof:}
First assume that the Lie algebra of $A\p$, hence $\Lie((A\p)\d)\d$,
satisfies the 
determinant condition, which means that $\Lie((A\p)\d)\d$ is an element of
$\N^{K}(S\p)$. 
By Theorem 2.9 it is an element of $\N^{DP}(S\p)$, i.e., $\Lie((A\p)\d)\d$ is
totally isotropic with respect to a perfect alternating $\Of$-linear form
on $\H(A\p/S\p)$. As all alternating $\Of$-linear forms on $\H(A\p/S\p)$
differ only by a scalar, the module $\Lie((A\p)\d)\d$ is totally
isotropic with respect 
to the alternating form induced by $\lambda$, hence $\lambda$ lifts to
$A\p$.
\hfill{$\Box$}

\begin{prop}
Let $A$ be an abelian scheme of type (DP). Then $\Lie(A)$ satisfies the
determinant condition.
\end{prop}

\noindent
\textsc{Proof:} Since $A$ is of type (DP), there exists a $p$-principal
polarization on $A$ locally for the \'etale topology on $S$. The determinant
condition can be proved locally for the \'etale topology, hence we may
assume that there 
exists an $\Of$-linear $p$-principal polarization $\lambda$ on $A$. As
in the case of the 
Tate module (proof of Proposition 1.3), the $p$-principal
polarization $\lambda$ induces a perfect alternating $\Os$-linear form on
$\H(A/S)$ such that the elements of $\Of$ are self-adjoint. By 1.5 we
obtain a perfect alternating $\Of\otimes_{\mathbb{Z}}\Os$-linear form
$\psi_{\lambda}$ on $\H(A/S)$.
Using the Hodge
filtration 
$$0\rightarrow\Lie(A\d)\d\rightarrow\H(A/S)\rightarrow\Lie(A)\rightarrow 0$$
it is easy to see that $\Lie(A\d)\d$ is totally isotropic with respect to
$\psi_{\lambda}$. The first De~Rham homology $\H(A/S)$ is locally on
$S$ a free $\Of\otimes_{\mathbb{Z}}\Os$-module of 
rank 2 ([14] Lemma 1.3). We may assume $\H(A/S)$ is free. In this case 
all $\Of\otimes_{\mathbb{Z}}\Os$-bilinear alternating
forms differ only by a scalar, hence $\Lie(A\d)\d$ is an element of the
local model $\N^{DP}(S)$. By Theorem 2.9 the local model $\N^{DP}$ is
equal to $\N^{K}$, 
thus $\Lie(A\d)\d$ and therefore $\Lie(A)$ satisfy the determinant condition. 
\hfill{$\Box$}\\

\noindent
\textbf{Remark 2.12.}
If $p$ is not equal to 2, it is easy to show that $\N^{DP}$ is contained in
$\N^{K}$, hence the determinant condition follows easily from condition
(DP). In case of $p=2$ we have to use the flatness of $\N^{DP}$ to show
this inclusion (cf. Section 5). 

\addtocounter{prop}{1}
\begin{prop}
Let $A$ be an abelian scheme which satisfies the determinant
condition. Then $A$ is of type (DP).
\end{prop}

\noindent
\textsc{Proof:} 
We will first give two different proofs of the proposition in the case of 
an abelian variety and then prove the general case.
Now let $A$ be an abelian variety over a field $k$.
Since all
conditions are trivial in characteristic zero,
we may assume
that the characteristic of $k$ is $p$. 
As
$\P(A)_{k}$ is equal to $\P(A)_{\overline{k}}$ (comp. Proposition 1.3
c)), we can assume that 
$k$ is an algebraically closed field. Furthermore, we will assume for
simplicity of notation that there is only one prime ideal
$\mathfrak{p}$ over $p$ in $\Of$.  

Yu has proved  the existence of an $\Of$-linear $p$-principal polarization
on $A$ by using Dieudonn\'e theory ([16] Lemma 2.6). We will sketch
the proof. Denote by 
$\widehat{\mathcal{O}}_{F}$ the completion of $\Of$ in
$\mathfrak{p}$. The Dieudonn\'e 
module $M$ of $A$ is a free
$\widehat{\mathcal{O}}_{F}\otimes_{\mathbb{Z}_{p}}W(k)$-module of 
rank 2 ([14] Lemma 1.3). Since $A$ is an abelian variety, there exists an
$\Of$-linear polarization $\lambda$ on $A$ ([14] Proposition 1.12). It
induces a non degenerate $\Of$-linear alternating form 
$$\psi_{\lambda}:M\times M\rightarrow D^{-1}\otimes_{\mathbb{Z}_{p}}W(k),$$
where $D^{-1}$ is the inverse different of $\widehat{F}_{\mathfrak{p}}$ over
$\mathbb{Q}_{p}$. 

If $p$ is totally ramified, denote by
$\pi$ a uniformizing element of $\widehat{\mathcal{O}}_{F}$. Then
$\widehat{\mathcal{O}}_{F}\otimes_{\mathbb{Z}_{p}}W(k)$ is a discrete
valuation ring with uniformizing element $\pi$ and $D^{-1}$ is equal
to $(\pi^{-d})$ for a positive integer $d$.  
Let $y_{1},y_{2}$ be an
$\widehat{\mathcal{O}}_{F}\otimes_{\mathbb{Z}_{p}}W(k)$-basis of $M$
and let $r$ be the order 
$\ord_{\pi}(\psi_{\lambda}(y_{1},y_{2}))$. Choose a totally positive
element $x\in F$ with $\ord_{\pi}(x)=-(r+d)$. Then the alternating
form  $x\psi_{\lambda}$ 
is a perfect integral form. Therefore,
the morphism $x\lambda$ of the abelian variety is well defined, hence
a $p$-principal 
polarization.

In the case of arbitrary ramification we obtain a decomposition of $M$
analogous to the 
decomposition of the Lie algebra (see 2.6). The determinant condition
implies that 
the order of $\psi_{\lambda}$ is equal on each component and we can choose
an element $x$ as above.

Now we will give an alternative proof of the proposition for an abelian
variety over an algebraically closed field $k$.
The idea is to show the existence of an $\Of$-linear $p$-principal
polarization of $A$ by lifting 
the abelian variety to characteristic zero (comp. [7] Lemma 5.5 and [16]
Example 8). Let $R$ be the ring
$W(k)[\sqrt{p}]$, where $W(k)$ is the ring of Witt vectors and let
$\mathfrak{m}$ be the maximal ideal of $R$. Let $R_{n}$ be the ring
$R/\mathfrak{m}^{n+1}$ for any positive integer $n$. Then $R_{n}$ is a
thickening of $R_{n-1}$ with trivially nilpotent PD-structure. 
Since $A$ is an abelian variety with $\Of$-operation $\iota$ over a
field $k$, there exists an  
$\Of$-linear polarization
$\lambda$ on $A$ ([14] Proposition 1.12). We will first lift
$(A,\iota,\lambda)$ to $R_{1}$. 

We have a decomposition
$$\Lie(A\d)\d=\bigoplus_{i=1}^{f}L_{i},$$
where $L_{i}$ is a $k[\pi]/(\pi^{e})$-submodule of the summand $H_{i}$ of
the decomposition $\H(A/k)=\bigoplus_{i=1}^{f}H_{i}$. There exists a
$k[\pi]/(\pi^{e})$-basis $f_{1}^{i},f_{2}^{i}$ of $H_{i}$ such that  
$$L_{i}=\langle\pi^{e_{1}^{i}}f_{1}^{i},\pi^{e_{2}^{i}}f_{2}^{i}
\rangle_{k[\pi]/(\pi^{e})},$$ 
and the determinant condition is equivalent to $e_{1}^{i}+e_{2}^{i}=e$
for every index $i$ (see 2.6). Hence $L_{i}$ is totally isotropic with
respect to every $k[\pi]/(\pi^{e})$-bilinear alternating form on
$H_{i}$ which is necessary for the existence of an $\Of$-linear
$p$-principal polarization.  

Let $\D(A/k)$
be the crystal associated to $A$. The module $\D(A/k)_{R_{1}}$ is a free
$\Of\otimes_{\mathbb{Z}}R_{1}$-module of rank 2 ([14] Lemma 1.3).
We have 
$$\Of\otimes_{\mathbb{Z}}W=\prod_{i=1}^{f}\mathcal{O}_{F,i},$$
where $\mathcal{O}_{F,i}$ is a totally ramified extension of $W$ of rank $e$. 
We obtain a decomposition of 
the crystal $\mathbb{D}(A/k)_{R_{1}}=\bigoplus_{i=1}^{f}D_{i}$, where
$D_{i}$ is a free $\mathcal{O}_{F,i}\otimes_{W}R_{1}$-module of rank 2. 

We want to use the theorem of Grothendieck-Messing to lift $(A,\iota,\lambda)$
to $R_{1}$. Therefore, we need to construct an
$\Of\otimes_{\mathbb{Z}}R_{1}$-module $\F$ of $\D(A/k)_{R_{1}}$, a
direct summand as an 
$R_{1}$-module, such that $\F$ lifts the Lie algebra and is totally
isotropic with respect to one, hence to every perfect alternating
$\Of\otimes_{\mathbb{Z}}R_{1}$-bilinear form on $\D(A/k)_{R_{1}}$. It is
sufficient to construct an $\mathcal{O}_{F,i}\otimes_{W}R_{1}$-module
$\F_{i}\subset D_{i}$ for every $i$ which lifts $L_{i}$ and which is
totally isotropic with respect to a perfect
$\mathcal{O}_{F,i}\otimes_{W}R_{1}$-bilinear 
form on $D_{i}$. For simplicity of notation we will write $\Of$
 and $\F$ instead of $\mathcal{O}_{F,i}$ and $\F_{i}$ etc.

Let
$f_{1},f_{2}$ be a $k[\pi]/(\pi^{e})$-basis of 
$H$ such that
$L=\langle\pi^{e_{1}}f_{1},\pi^{e_{2}}f_{2}\rangle_{k[\pi]/(\pi^{e})}$
(see 2.6). 
Since $\mathcal{O}_{F}$ is a totally ramified extension of $W$,
we 
have
$\mathcal{O}_{F}=W[\pi]/(P(\pi)),$
 where $P$ is an Eisenstein polynomial of degree $e$. 
Let $P(\pi)$ be
$\pi^{e}-px$ with $x\in\Of^{\times}$. Define 
$$\mathcal{F}=\langle\pi^{e_{1}}f_{1}+x\sqrt{p}f_{2},\pi^{e_{2}}f_{2}+
\sqrt{p}f_{1}\rangle_{\mathcal{O}_{F}\otimes_{W}R_{1}}.$$
The submodule $\mathcal{F}$ is $\pi$-invariant and lifts
$L$. Furthermore, $\mathcal{F}$ is   
totally isotropic with
respect to the $\mathcal{O}_{F}\otimes_{W}R_{1}$-linear
alternating form $\langle\cdot,\cdot\rangle$ defined by $\langle
f_{1},f_{2}\rangle=1$ because the determinant condition implies that
$e_{1}+e_{2}$ is equal to $e$.

Using the theorem of Grothendieck-Messing ([9] 4.1.10),
we obtain an abelian scheme $\tilde{A}_{1}$ over $R_{1}$ with
$\Of$-operation $\tilde{\iota}_{1}$ and $\Of$-linear
polarization $\tilde{\lambda}_{1}$. By induction we can lift
$(\tilde{A}_{n},\tilde{\iota}_{n},\tilde{\lambda}_{n})$ to
$(\tilde{A}_{n+1},\tilde{\iota}_{n+1},\tilde{\lambda}_{n+1})$ over
$R_{n}$ using the same 
arguments as above.
We obtain a compatible
system of liftings $(\tilde{A}_{n},\tilde{\iota}_{n},\tilde{\lambda}_{n})$
of $(A,\iota,\lambda)$ 
over $R_{n}$ for every integer $n$.
The polarizations $\tilde{\lambda}_{n}$ induce a compatible system of ample
invertible sheaves on every abelian scheme $\tilde{A}_{n}$.
Therefore, the projective limit of all
$\tilde{A}_{n}$ is algebraized by an abelian scheme $\tilde{A}$
with $\Of$-operation $\tilde{\iota}$ and $\Of$-linear polarization
$\tilde{\lambda}$ over the discrete valuation ring $R$. The generic
fibre $\tilde{A}_{\eta}$ is an 
abelian variety over a field of characteristic zero. Hence there
exists an $\Of$-linear $p$-principal polarization $\mu_{\eta}$ on
$\tilde{A}_{\eta}$. We have  
$$\Hom(\tilde{A}_{\eta},\tilde{A}_{\eta}\d)=\Hom(\tilde{A},\tilde{A}\d)$$
because $\tilde{A}\d$ is the N\'eron model of $\tilde{A}_{\eta}\d$ ([3]
1.4.2). Thus $\mu_{\eta}$ extends to an $\Of$-linear isogeny $\mu$ prime to
$p$. Since $\tilde{\lambda}_{\eta}$ and $\mu_{\eta}$ differ by a totally
positive element of $F$, the isogeny $\mu$ is a polarization. 
After reduction modulo $p$ we obtain an $\Of$-linear $p$-principal
polarization of $A$.\\

Now we will prove the proposition in the general case of an abelian scheme $A$.
We may assume that $A$ is defined over a
local ring $R$ with residue field $k$ of characteristic $p$ such that $R$
is the localization of a finitely generated $\mathbb{Z}$-algebra. We have
to show that there exists an $\Of$-linear $p$-principal polarization
over the henselization 
of $R$. Using the approximation theorem of Artin for $\P(A)$ ([1]
Theorem 1.12), we only have 
to prove the existence over the completion $\hat{R}$ of $R$. Hence we may
assume that $R$ is a complete local ring. Let $k$ be the residue field of $R$. 
We have already proved
that there exists an $\Of$-linear $p$-principal polarization $\lambda$
on the abelian variety 
$A_{k}$. 

We will construct an $\Of$-linear $p$-principal polarization on $R$ by lifting
$\lambda$ inductively to an $\Of$-linear  $p$-principal polarization
$\lambda_{n}$ over 
$R/\mathfrak{m}^{n}$. Then the
compatible system of $\lambda_{n}$ over
$R/\mathfrak{m}^{n}$ will induce an $\Of$-linear $p$-principal
polarization on $A$ by Grothendieck's algebraization theorem. Consider
the surjection $R/\mathfrak{m}^{n+1}\rightarrow
R/\mathfrak{m}^{n}$. The square of the  ideal
$(\mathfrak{m}^{n+1}/\mathfrak{m}^n)$ in $R/\mathfrak{m}^{n+1}$ is
equal to zero, hence $(\mathfrak{m}^{n+1}/\mathfrak{m}^n)$ can be
endowed with the trivial PD-structure. By Corollary 2.10 we can lift
$\lambda_n$ to an $\Of$-linear $p$-principal polarization
$\lambda_{n+1}$ of $A_{R/\mathfrak{m}^{n+1}}$ which proves the
proposition.  
\hfill{$\Box$}\\

\begin{prop}
Let $(A,\iota)$ be an abelian variety of type (DP) over
a field $k$ of characteristic $p$ and let $(\tilde{A},\tilde{\iota})$
be a deformation of $(A,\iota)$ over a local artinian ring $R$.  
Then $\P(\tilde{A})$ is locally constant if and only if $\tilde{A}$ is
of type (DP). 
\end{prop}

\noindent
\textsc{Proof:}
If $\P(\tilde{A})$ is locally constant, the canonical morphism
$$\tau:\tilde{A}\otimes_{\Of}\P(\tilde{A})\rightarrow \tilde{A}\d$$
is a morphism of abelian schemes (Lemma 1.9). Since $A$ is of type
(DP), the fibre morphism $\tau_k$, and hence $\tau$, is an
isomorphism.  

The converse implication is true by Proposition 1.4.
\hfill{$\Box$}\\
  
\noindent
\textbf{Remark 2.15.}
Let $(A,\iota)$ be an abelian variety of type (DP) over
a field $k$ of characteristic $p$.
Denote by $\Def[A,\iota]$  the equi-characteristic deformation functor
of $(A,\iota)$ and denote by $\Def[A,\iota,\lambda]$  the analogous
deformation functor of the tuple $(A,\iota,\lambda)$ where $\lambda$
is an $\Of$-linear $p$-principal polarization of $A$.  
Yu has proved ([16]
Remark 2.18(2)) that the morphism
$$\Def[A,\iota,\lambda]\rightarrow\Def[A,\iota]$$
is formally \'etale if and only if
$A$ is of type (R). Hence for every abelian variety $A$ of type (DP)
and not of type (R) there exists a deformation $\tilde{A}$ not of type
(DP), i.e., where $\P(\tilde{A})$ is not locally constant (Proposition
2.14). 
\\

\noindent
\textbf{Example 2.16.}
We now give an example of an abelian scheme, where the polarization
module $\P$ is not locally 
constant. 
Let $[F:\mathbb{Q}]$ be equal to 2 and let
$p$ be totally ramified in $\Of$.
Let $A$ be
an abelian variety not of type (R) over a field $k$ of characteristic
$p$. Such an abelian variety exists (see Remark 3.8). 
We have $\Of\otimes_{\mathbb{Z}}k=k[\pi]/(\pi^{2})$. Since $A$ is not of
type (R), there exists a
basis $f_{1}, f_{2}$ of the free $k[\pi]/(\pi^{2})$-module $\H(A/k)$ such that
$\Lie(A\d)\d$ is generated by $\pi f_{1}$ and $\pi f_{2}$ as a
$k[\pi]/(\pi^{2})$-module (see 2.6). Since the determinant condition
in the totally 
ramified case is trivial for an abelian variety over a field of characteristic
$p$ (cf. 2.6), there exists an $\Of$-linear $p$-principal polarization
on $A$ (Proposition 2.13).  

Let $\D(A/k)$ be the crystal to $A$.
We want to construct
a lift $\tilde{A}$ of $A$ over $k[\epsilon]=k[X]/(X^{2})$ such that
$\Lie(\tilde{A}\d)\d$ is not totally isotropic with respect to a
perfect alternating 
$\Of\otimes_{\mathbb{Z}}k[\epsilon]$-linear form $\langle\cdot,\cdot\rangle$ on
$\D(A/k)_{k[\epsilon]}$. Every $\Of$-linear $p$-principal polarization
induces a 
perfect alternating
$\Of\otimes_{\mathbb{Z}}k[\epsilon]$-linear form on
$\D(A/k)_{k[\epsilon]}$ which differs from $\langle\cdot,\cdot\rangle$
only by a 
unit. Thus there cannot exist an $\Of$-linear $p$-principal polarization on
$\tilde{A}$. This means that $\P(\tilde{A})$ is not constant and
$\tilde{A}\otimes_{\Of}\P(\tilde{A})$ is not represented by an abelian
scheme. Indeed, otherwise it would follow that the isogeny 
$$\tau:\tilde{A}\otimes_{\Of}\P(\tilde{A})\rightarrow\tilde{A}\d$$
is an isomorphism, i.e., $\tilde{A}$ is of type (DP), which is impossible since
there exists no $\Of$-linear $p$-principal polarization on
$\tilde{A}$. 

To construct $\tilde{A}$ denote again by $f_{1}, f_{2}$ a lift of the
$k[\pi]/(\pi^{2})$-basis $f_{1}, f_{2}$ of $\H(A/k)$ to a
$k[\epsilon,\pi]/\pi^{2}$-basis of 
$\D(A/k)_{k[\epsilon]}$. Let $\mathcal{F}$ be the module
$$\mathcal{F}=\langle (\pi+\epsilon) f_{1},\pi f_{2}\rangle_{k[\epsilon]}.$$
It is a free $k[\epsilon]$-module of rank 2 which is $\pi$-invariant and
lifts $\Lie(A\d)\d$. Let $\langle\cdot,\cdot\rangle$ be the canonical
alternating 
$\Of\otimes_{\mathbb{Z}}k[\epsilon]$-linear form associated to the basis
$f_{1}, f_{2}$. We obtain $\langle(\pi+\epsilon)f_{1},\pi
f_{2}\rangle=\pi\epsilon$, hence $\F$ is not totally isotropic with respect
to $\langle\cdot,\cdot\rangle$. Using the theorem of Grothendieck-Messing we
obtain the desired
abelian scheme over $k[\epsilon ]$.


\section{The moduli problems of Deligne/Pappas and Kottwitz}
In this section we will introduce two different integral moduli
problems for the 
Shimura variety corresponding to the group $G$ as in the
introduction. We will show in 
Section 4 that these two moduli problems are isomorphic. 

Let $(L,L^{+})$ be an invertible $\Of$-module with
positivity and let $(A,\iota)$ be an abelian scheme with $\Of$-operation over a
base scheme $S$. 

\begin{de} (Deligne/Pappas)
An $L$-polarization of $(A,\iota)$ is a homomorphism of $\Of$-modules
\begin{align}
\varphi:L&\rightarrow\P(A)_{S}\notag\\
\lambda&\mapsto\varphi_{\lambda},\notag
\end{align}
which is compatible with the positivities on $L$ and $\P(A)_{S}$
respectively, such that 
the composition
$$A\otimes_{\Of}L\stackrel{1\otimes\varphi}{\longrightarrow}
A\otimes_{\Of}\P(A)_{S}\stackrel{\tau}{\longrightarrow} A\d$$ 
is an isomorphism.
\end{de}

\noindent
\textbf{Remark 3.2.} Let $A$ be an abelian scheme with $\Of$-operation
and $L$-polarization $\varphi$. 
Consider the exact sequence (1.8.1)
$$0\longrightarrow A[\mathfrak{a}_{\lambda}]\longrightarrow A
\stackrel{\lambda}{\longrightarrow}A\otimes_{\Of}L\longrightarrow 0,$$ 
where $\mathfrak{a}_{\lambda}$ is the annihilator of the $\Of$-module
$L/\lambda\Of$. It is clear that the morphism $A\otimes_{\Of}L\rightarrow
A\d$ is an isomorphism if and only if for one and thus for all elements
$\lambda\neq 0$ of $L$ the morphism $\varphi_{\lambda}$ is an isogeny with
kernel equal to $A[\mathfrak{a}_{\lambda}]$. 

\addtocounter{prop}{1} 
\begin{prop}
A homomorphism $\varphi:L\rightarrow\P(A)_{S}$ of $\Of$-modules is an
$L$-polarization of $(A,\iota)$ if and only if $A$ is of type (DP) and the
induced morphism
$\varphi:(L,L_{+})\rightarrow(P(A),\P(A)_{+})$ of sheaves for the \'etale
topology is an isomorphism.
\end{prop}
 
\noindent
\textsc{Proof:}
Suppose that $A$ is of type (DP) and $\varphi$ is an isomorphism. We obtain 
$$A\otimes_{\Of}L\stackrel{\sim}{\longrightarrow}A\otimes_{\Of}\P(A)
\stackrel{\sim}{\longrightarrow}A\d,$$ 
hence $\varphi$ is an $L$-polarization.

Conversely, suppose that $\varphi$ is an $L$-polarization. It induces a
morphism of sheaves $\varphi:L\rightarrow\P(A)$. 
Because the bijectivity of
$\varphi$ is a local property and $S$ is locally noetherian, we may
assume that $S$ is connected. It is 
sufficient to show that $\varphi:L\rightarrow P:=\P(A)_{S}$ is an isomorphism. 
It is clear that $\varphi$ is injective.
We can prove the surjectivity of $\varphi$ locally on $\Of$,
hence we can assume that $P$ and $L$ are free $\Of$-modules. Let
$\lambda$ be a generator of $P$ and let $x\lambda$ be a generator of $L$,
where $x$ is an element of $\Of$. 
The isomorphism $A\otimes_{\Of}L\stackrel{\sim}{\longrightarrow}A\d$
factors through $A\otimes_{\Of}P$. Hence we obtain an injective morphism
$$A\otimes_{\Of}(x\lambda)\hookrightarrow A\otimes_{\Of}(\lambda).$$
Thus the kernel $A[x]$ is equal to zero. But the degree of the
multiplication with $x$ on $A$ is equal to $N(x)^{2}$. Therefore, $x$
is an element of $\Of^{\times}$ and $\varphi$ is an 
isomorphism.

This implies that the
isogeny $\tau: A\otimes_{\Of}\P(A)\rightarrow A\d$ is an isomorphism and that
$\P(A)$ is constant, hence
$A$ is of type (DP). \hfill{$\Box$}\\

\noindent
\textbf{Remark 3.4.}
Rapoport considered in [14] abelian schemes of type (R) together with
an isomorphism 
$\varphi:(L,L_{+})\rightarrow(P(A),\P(A)_{+})$ of sheaves for the \'etale
topology  of $\Of$-modules. Since every abelian scheme
of type (R) is of type (DP),
the last proposition shows that it is the same
to consider abelian schemes of type (R) with $L$-polarization.\\

\noindent
We will now concentrate on the case where $L$ is equal to the inverse
different $D^{-1}$ with positivity induced from the totally 
positive elements in $F$.

\addtocounter{prop}{1}
\begin{de}
Let $A$ be an abelian scheme with $\Of$-operation $\iota$ and
$D^{-1}$-polarization $\varphi$ and let $n$ be an integer prime to $p$. 
A level $n$ structure of $A$ is an $\Of$-linear isomorphism of
commutative group schemes over $S$  
$$\tau :A[n]\stackrel{\sim}{\longrightarrow}(\Of/n\Of)^{2}$$
which satisfies the following
condition. There exists an
isomorphism
$f:\mu_{n}\stackrel{\sim}{\longrightarrow}\mathbb{Z}/n\mathbb{Z}$
such that the diagram
$$\xymatrix{
\P(A)\otimes_{\Of}\wedge^{2}_{\Of/n\Of}A[n]
\ar@{->}[r]_-{e_{n}-pairing}^-{\sim}\ar@{->}[d]^{\sim}_{\varphi^{-1}
  \otimes\tau}  
& D^{-1}\otimes_{\mathbb{Z}}\mu_{n}\ar@{->}[d]^{\sim}_f\\ 
D^{-1}\otimes_{\Of}\wedge^{2}_{\Of/n\Of}(\Of/n\Of)^{2} \ar@{->}[r]^-{\sim}
& D^{-1}\otimes_{\mathbb{Z}}\mathbb{Z}/n\mathbb{Z}\\
}$$
commutes.
\end{de}
 
\noindent
\textbf{3.6. The moduli problem of Deligne and Pappas}\\
Fix an integer $n\geq 3$ prime to $p$ and 
consider the moduli problem $\Md$ over $\Spec\mathbb{Z}_{(p)}$ which
associates to every scheme $S$ over $\Spec\mathbb{Z}_{(p)}$ the following
data up to isomorphism
\begin{enumerate}
\item An abelian scheme $A$ over $S$ of relative dimension $g$.
\item An $\Of$-operation $\iota:\Of\rightarrow \End(A)$ on $A$.
\item An $D^{-1}$-polarization $\varphi:D^{-1}\rightarrow\P(A)_{S}$. 
\item A level $n$ structure
  $\tau:A[n]\stackrel{\sim}{\longrightarrow}(\Of/(n\Of))^{2}$. 
\end{enumerate}
\textbf{Remark 3.7.} 
The functor $\Md$ is represented by a separated scheme of finite type over
$\Spec\mathbb{Z}_{(p)}$. \\

\noindent
\textbf{Remark 3.8.}
Rapoport considered in [14] the subscheme $\Mr$ of all abelian schemes
in $\Md$ of type (R). This subscheme is the smooth locus of $\Md$ 
([14] Theorem 1.20, [16] Proposition 2.15).

Since condition (R) is trivial in characteristic zero, the generic fibre of
$\Md$ is equal to the generic fibre of $\Mr$. The moduli problems are the
same if $p$ is unramified because in this case conditions (R) and (DP) are
equivalent. If $p$ ramifies, the singular locus of $\M^{DP}$ is of
codimension 2 ([5]
Theorem 2.2), hence
there exist abelian varieties of type (DP) over an algebraically closed
field of characteristic $p$ such that the Lie algebra is not free. \\

\noindent
\textbf{3.9.}
Another integral model for the Shimura variety of the group $G$ as in the
introduction is the specialization to
the Hilbert-Blumenthal case of the moduli problem of
Kottwitz type defined by Kottwitz [8] in the unramified case and Rapoport and
Zink [15] in the general case.

We fix a totally positive element $d\in D^{-1}$ such that
$\mid\! D^{-1}/d\Of\!\mid$ is prime to $p$.
Denote by $\Zhp$ the product $\prod_{l\neq p}\mathbb{Z}_{l}$
and by $\Ap$ the ring $\Zhp\t \mathbb{Q}$. Furthermore, let 
$$\Tp=\prod_{l\neq p}T_{l}(A)$$ 
be
the product of all $l$-adic Tate modules.
Denote by $\Gn$ the principal congruence subgroup modulo $n$
$$\Gn=\{g\in Gl_{2}(\Of\t\Zhp)\mid \det(g)\in(\Zhp)^{\times}\text{ and }
g\equiv 1 \text{ mod } n\}$$
and let $\langle\cdot,\cdot\rangle$ be the canonical
$\Of\otimes_{\mathbb{Z}}\Zhp$-bilinear alternating form on
$(\Of\otimes_{\mathbb{Z}}\Zhp)^{2}$.
We want to define a level structure of
type $\Gn$
for an abelian scheme up to isogeny prime to $p$ with
$\Of$-operation and $\Of$-linear $p$-principal polarization. It is
sufficient to do this in the case of a connected base 
scheme $S$. 

\addtocounter{prop}{4}
\begin{de}
(Kottwitz) Let $A\otimes_{\mathbb{Z}}\mathbb{Z}_{(p)}$ be an abelian
scheme up to 
isogeny prime to $p$ with
$\Of$-operation and $\Of$-linear $p$-principal
polarization $\lambda$ over a connected scheme $S$. Let $s$ be a geometric
point of $S$. The group $\Gn$ operates on the set of $\Of\t\Ap$-linear
isomorphisms 
$$\widehat{T}^{p}(A_{s})\t\mathbb{Q}\stackrel{\sim}{\longrightarrow}
(\Of\t\Ap)^{2}.$$  
A level structure of type 
$\Gn$ is a class modulo $\Gn$ of $\Of\t\Ap$-linear isomorphisms
$\overline{\eta}^{p}$ such that the alternating form induced by $\lambda$ on
the left hand side is equal to the form $d\langle\cdot,\cdot\rangle$
on the right 
hand side up to an element of $(\Ap)^\times$. Furthermore, the class
of isomorphisms 
has to be invariant under the fundamental group. Here $d$
is the fixed element of $D^{-1}$.
\end{de}

\noindent
\textbf{3.11. The moduli problem of Kottwitz type}\\
Fix an integer $n\geq 3$
prime to $p$ and
consider the moduli problem $\Mk$ over $\Spec\mathbb{Z}_{(p)}$
which associates to every scheme $S$ over $\Spec\mathbb{Z}_{(p)}$ the
following data up to isomorphism
\begin{enumerate}
\item An abelian scheme $A\t\mathbb{Z}_{(p)}$ over $S$ of relative
  dimension $g$ up to isogeny prime to $p$ such that the Lie algebra
  of $A\otimes_{\mathbb{Z}}\mathbb{Z}_{(p)}$ 
satisfies the determinant condition.
\item A homomorphism of rings $\iota\t\mathbb{Z}_{(p)}:\Of\rightarrow
  \End(A)\t\mathbb{Z}_{(p)}$. 
\item A $\mathbb{Q}$-subspace $\overline{\lambda}$ of
  $(\P(A)_{S}\t\mathbb{Z}_{(p)})\otimes_{\mathbb{Z}_{(p)}}\mathbb{Q}$ of
  dimension one which contains a $p$-principal polarization.
\item A level structure of type $\Gn$, 
$$\overline{\eta}^{p}:\Tp\t\mathbb{Q}\stackrel{\sim}{\longrightarrow}
(\Of\t\mathbb{A}_{f}^{p})^{2}\mod\Gn.$$ 
\end{enumerate}

\noindent
\textbf{Remark 3.12.} 
The moduli problem $\Mk$ is represented by a separated scheme of finite
type over $\Spec\mathbb{Z}_{(p)}$ ([8] Section 5).

The moduli problem $\Mk=\Mk_{d}$ does not depend on the choice of the
totally positive 
element $d\in D^{-1}$.
Let $d\p$ be another totally positive element of $D^{-1}$ such that $\mid\!
D^{-1}/d\Of\!\mid$ is prime to $p$. Then $d\p$ is equal to $xd$, where $x$
is a totally positive element of $F$, and we obtain a canonical isomorphism
\begin{align}
\sigma:\Mk_{d}&\stackrel{\sim}{\longrightarrow}\Mk_{d\p}\notag\\
(A\t\mathbb{Z}_{(p)},\iota\t\mathbb{Z}_{(p)},\overline{\lambda},
\overline{\eta}^{p})&\mapsto(A\t\mathbb{Z}_{(p)},\iota\t\mathbb{Z}_{(p)},
\overline{x\lambda},\overline{\eta}^{p}).\notag 
\end{align}
Thus $\Mk_{d}$ is independent of $d$ up to canonical
isomorphism and we will write $\Mk$ instead of $\Mk_{d}$.

The determinant condition is an automatic consequence of
the existence of an $\Of$-linear $p$-principal polarization. Indeed, the Lie
algebra is totally isotropic with respect to the alternating form induced
by the $p$-principal polarization, hence it is an element of the local
model $\N^{DP}$ (see 2.8). Since $\N^{DP}$ is equal to the local model
$\N^{K}$ defined by the determinant condition (Theorem
2.9),
 the Lie algebra satisfies automatically the determinant condition.


\section{Comparison of the two moduli problems}

In this section we will construct an isomorphism between the moduli
problems $\Md$ and $\Mk$.

Let $n$ be an integer prime to $p$ and let $d\in D^{-1}$ be a
totally positive element such that $\mid\! D^{-1}/d\Of\!\mid$ is prime
to $p$. Denote by $\langle\cdot,\cdot\rangle$
the canonical alternating $\Of\t\Zhp$-bilinear form on $(\Of\t\Zhp)^{2}$.

\begin{lem}
Let $(A,\iota,\varphi,\tau)$ be an element of $\Md(S)$,
where $S$ is connected. Then there exists an isomorphism 
$$\eta^{p}: \Tp\stackrel{\sim}{\longrightarrow}(\Of\t\Zhp)^{2}$$
which lifts $\tau$ such that the alternating forms $\psi_{\varphi(d)}$
on the Tate 
module and $d\langle\cdot,\cdot\rangle$ on
the right hand side only differ by an
element of $(\Zhp)^{\times}$. The isomorphism is uniquely determined up to
an element of $\Gn$. In particular, $\tau$ induces a uniquely determined
level structure of type $\Gn$ for
$(A\otimes_{\mathbb{Z}}\mathbb{Z}_{(p)},\iota\otimes_{\mathbb{Z}}
\mathbb{Z}_{(p)},\varphi(d)\mathbb{Q})$.   
\end{lem}

\noindent
\textsc{Proof:}
First suppose that $S$ is the spectrum of an algebraically closed field. 
The isomorphism 
$$\tau:A[n]\stackrel{\sim}{\longrightarrow}(\Of/n\Of)^{2}$$ 
is
a level $n$ structure, thus the restricted alternating forms
$\psi_{\varphi(d),n}$ and $d\langle\cdot,\cdot\rangle_{n}$ only differ
by an element of $(\mathbb{Z}/n\mathbb{Z})^{\times}$. Therefore, we
can lift $\tau$ to an 
isomorphism 
\begin{align}
\eta^{p}: \Tp\stackrel{\sim}{\longrightarrow}(\Of\t\Zhp)^{2}\tag{4.1.1}
\end{align}
such that the corresponding forms $\psi_{\varphi(d)}$ and
$d\langle\cdot,\cdot\rangle$ on 
the Tate module only differ by an
element of $(\Zhp)^{\times}$. The
lift is uniquely determined up to an element of $\Gn$. 

In the general case let $s$ be a geometric point of $S$. The fundamental group
acts trivially on $A_{s}[n]$. Thus the operation of an element of
$\pi_{1}(S,s)$ on the left hand side of (4.1.1) corresponds to the
operation of an 
element of $\Gn$ on the right hand side. Hence $\overline{\eta}^{p}$ is
invariant under $\pi_{1}(S,s)$. 
\hfill{$\Box$}\\

\noindent
\textbf{4.2.}
Using Lemma 4.1 we obtain a morphism
\begin{align}
\Phi:\Md&\rightarrow\ \Mk\notag\\
(A,\iota,\varphi,\tau)&\mapsto(A\t\mathbb{Z}_{(p)},\iota\otimes
\mathbb{Z}_{(p)},\varphi(d)\mathbb{Q},\overline{\eta}^{p}),\notag 
\end{align}
where $\overline{\eta}^{p}$ is the uniquely determined level structure of
type $\Gn$ which lifts $\tau$. Since $A$ is of type (DP), it
satisfies the Kottwitz determinant condition (Proposition 2.11). 

If we choose
another element $d\p\in D^{-1}$, we obtain a morphism 
$$\Phi_{d\p}:\Md\rightarrow\Mk_{d\p}$$
which is equal to $\sigma\circ\Phi$, where
$\sigma:\Mk_{d}\stackrel{\sim}\longrightarrow\Mk_{d\p}$ is the
canonical isomorphism 
(Remark 3.12). Thus $\Phi$ is independent of $d$ up to the canonical
isomorphisms of the 
moduli problem $\Mk$.  
 
\addtocounter{prop}{1}
\begin{thm}
The morphism 
$$\Phi:\Md\rightarrow\Mk$$
is an isomorphism.
\end{thm}

\noindent
\textsc{Proof:}
It is sufficient to construct the inverse functor on every connected scheme
$S$.
Let
$X=(A\t\mathbb{Z}_{(p)},\iota\t\mathbb{Z}_{(p)},\overline{\lambda},
\overline{\eta}^{p})$ 
be an element of $\Mk(S)$ for a connected scheme $S$. 
Denote by $s$ a geometric
point of $S$. Consider the level structure
$$\overline{\eta}^{p}:
\widehat{T}^{p}(A_{s})\otimes_{\mathbb{Z}}\mathbb{Q}
\stackrel{\sim}{\longrightarrow}(\Of\t\Ap)^{2}\  
\text{mod}\ \Gn\text{.}$$
The lattice 
$\Lambda=(\overline{\eta}^{p})^{-1}((\Of\t\Zhp)^{2})\subset
\widehat{T}^{p}(A_{s})\otimes_{\mathbb{Z}}\mathbb{Q}$  
is independent of the representative of $\overline{\eta}^{p}$, in
particular it is invariant under
$\pi_{1}(S,s)$. Therefore, it induces an abelian scheme $B$ over $S$ with
Tate module $\Lambda$
([4] 4.2) which is isomorphic to $A$ up to isogeny prime to $p$. The
abelian scheme $B$ is uniquely determined up to isomorphism. 

Since the Tate
module of $B$ is $\Of$-invariant, we obtain an $\Of$-operation
$\iota:\Of\rightarrow\End(B)$. Furthermore, an $\Of$-linear
$p$-principal polarization 
$\lambda\in \overline{\lambda}$ on $A$ induces an $\Of$-linear $p$-principal
polarization $\mu$ on $B$ uniquely determined up to a positive element of
$\mathbb{Z}_{(p)}$. By  construction of $B$ the level structure
$\overline{\eta}^{p}$ induces an isomorphism
$$\widehat{T}^{p}(B_{s})\stackrel{\sim}{\longrightarrow}(\Of\t\Zhp)^{2}\
\text{mod}\ \Gn$$
such that the
alternating forms $\psi_{\mu}$ and $d\langle\ ,\
\rangle$ differ only by an element of $(\Ap)^{\times}$. 
Modulo $n$ we obtain an isomorphism
$$\tau:B[n]\stackrel{\sim}{\longrightarrow}(\Of/n\Of)^{2}\text{.}$$

To get an element of $\Md(S)$ we have to
construct a $D^{-1}$-polarization and verify the compatibility condition
with the level $n$ structure.
First note that the
alternating forms $\psi_{\mu}$ and $d\langle\ ,\
\rangle$ are perfect forms at almost all places. Therefore,
we obtain that
$\psi_{\mu}$ is equal to $(qu)d\langle\cdot,\cdot\rangle$, where $q$
is an element of 
$\mathbb{Q}$ and $u$ an element of $(\Zhp)^{\times}$. Thus there exists a
uniquely determined $\Of$-linear $p$-principal polarization
$q^{-1}\mu\in \mu\mathbb{Q}$ 
of $B$ such that the
induced bilinear form is equal to $d\langle\cdot,\cdot\rangle$ up
to an element of $(\Zhp)^{\times}$. We will denote the form $q^{-1}\mu$
again by $\mu$.
Consider the isomorphism of $F$-vector
spaces  
\begin{align}
\varphi:D^{-1}\otimes_{\Of}F&\stackrel{\sim}{\longrightarrow}\P(B)_{S}
\otimes_{\Of}F\notag\\ 
d&\longmapsto \mu\notag
\end{align}
defined on the bases $d$ and $\mu$ respectively.\\

\emph{
Claim:} $\varphi$ induces an isomorphism of $D^{-1}$ and $\P(B)_{S}$. \\

\noindent
Since there exists an $\Of$-linear $p$-principal polarization on $B$,
the polarization 
module $\P(B)$ is constant (Proposition 1.4) and $B$ is of type (DP),
hence we obtain from the claim a $D^{-1}$-polarization on $B$.\\

We now prove the claim.
Since
$\overline{\eta}^{p}$ is a level structure, the induced forms of
$\psi_{\lambda}$ and $d\langle\cdot,\cdot\rangle$ on the Tate module differ
only by an element of 
$(\Zhp)^{\times}$. Thus the morphism
\begin{align}
D^{-1}\t\Zhp&\rightarrow\P(B)_{S}\t\Zhp\notag\\
d&\mapsto\mu\notag
\end{align}
is an isomorphism.

On the other hand, the elements $d$ and $\mu$ generate
the $\Of\t\mathbb{Z}_{p}$-modules $D^{-1}\t\mathbb{Z}_{p}$ and
$\P(B)_{S}\t\mathbb{Z}_{p}$ respectively because
$\mid\!D^{-1}/d\Of\!\mid$ is prime to $p$ and $\mu$ is an $\Of$-linear
$p$-principal polarization. Thus $\varphi$ factors through both  
modules and is an isomorphism.
Furthermore, $d$ and $\mu$ are totally positive, thus $\varphi$ induces a
bijection between the sets of totally positive elements ([17] Lemma
1.4). Thus the claim is proved. 

We now want to prove that
the above isomorphism 
$$\tau:B[n]\stackrel{\sim}{\longrightarrow}(\Of/n\Of)^{2}$$
is a level $n$ structure of $B$.
By construction of $\varphi$ the element $d$ maps to the
polarization $\mu$ and the induced alternating forms $\psi_{\mu}$ and
$d\langle\cdot,\cdot\rangle$ on the Tate module differ only by an element of
$(\mathbb{Z}_{p})^{\times}$. Hence these forms differ by an element of
$(\mathbb{Z}/n\mathbb{Z})^{\times}$ modulo $n$ and we can find an
isomorphism $\mu_{n}\stackrel{\sim}{\longrightarrow}\mathbb{Z}/n\mathbb{Z}$
such that the diagram 
$$\xymatrix{
\P(B)_{S}\otimes_{\Of}\wedge^{2}_{\Of/n\Of}B[n]
\ar@{->}[r]_-{e_{n}-pairing}^-{\sim}\ar@{->}[d]^{\sim}_{\varphi^{-1}\otimes\tau} & D^{-1}\otimes_{\mathbb{Z}}\mu_{n}\ar@{->}[d]^{\sim}_f\\
D^{-1}\otimes_{\Of}\wedge^{2}_{\Of/n\Of}(\Of/n\Of)^{2}\ar@{->}[r]^-{\sim}
& D^{-1}\otimes_{\mathbb{Z}}\mathbb{Z}/n\mathbb{Z}\\
}$$
commutes.

Altogether we have defined a morphism
\begin{align}
\Psi:\ \Mk&\rightarrow\Md\notag\\
(A\t\mathbb{Z}_{(p)},\iota\t\mathbb{Z}_{(p)},\overline{\lambda},
\overline{\eta}^{p})&\mapsto(B,\iota,\varphi,\tau)\text{.}\notag 
\end{align}
It is easy to verify that $\Psi$ is the inverse morphism of $\Phi$.
\hfill{$\Box$}\\


\section{Comparison of the corresponding local models }

In this section we will prove that the local models
$\N^{DP}$ and $\N^{K}$ (cf. 2.8) are equal as
subschemes of the Grassmannian, which was used already in Section 2.
The proof is based on the flatness of $\N^{DP}$ and $\N^{K}$. In case of
$p\neq 2$ the inclusion $\N^{DP}\subset\N^{K}$ is elementary and gives a
new proof of the flatness of $\N^{DP}$ by using the flatness of
$\N^{K}$. Later we will give a direct proof of this equality in case of tame
ramification. \\

\noindent
\textbf{5.1.}
We first explain the relation of the moduli problems $\M^{DP}$ and $\M^{K}$
to the corresponding local models.
Let $(A,\iota,\varphi,\eta)$ be an element of $\M^{DP}(S)$. For an
$\Os$-module $M$ denote by
$M\d$ its dual module.
A trivialization
of $\H(A/S)$ is an isomorphism 
$$\Phi:(\Of\otimes_{\mathbb{Z}}\Os)^{2}\stackrel{\sim}{\longrightarrow}
\H(A/S)$$
such that the diagram
$$\xymatrix{
(\Of\otimes_{\mathbb{Z}}\Os)^{2}\otimes_{\Of}D^{-1}\ar@{->}[r]_-\sim^-{\Phi
  \otimes\varphi}\ar@{->}[d]^\sim_{tr_{F/\mathbb{Q}}}&\H(A/S)\otimes_{\Of} 
\P(A)\ar@{->}[d]^\sim_\tau\\    
((\Of\otimes_{\mathbb{Z}}\Os)^{2})\d\ar@{->}[r]_-\sim^-{\Phi\d}&\H(A\d/S)\\
}$$
commutes, where $\tau$ is the canonical isomorphism
$A\otimes_{\Of}\P(A)\stackrel{\sim}{\longrightarrow} A\d$ (Proposition
2.2).
Consider the functor $\Ms^{DP}$ of $\M^{DP}$ over $\Spec\mathbb{Z}_{(p)}$ with
\begin{align}
\Ms^{DP}(S)=\{(A,\iota,\varphi,\eta,\Phi)\mid
(A,\iota,\varphi,\eta)\in\M^{DP}(S)\text{ and }\notag\\
\Phi
\text{ a trivialization
  of }\H(A/S)\}.\notag
\end{align}
Let $\pi:\Ms^{DP}\rightarrow\M^{DP}$ be the natural forgetful functor and let
$G$ be the algebraic group over
$\mathbb{Z}$ given by
$$G(R)=\{g\in
\End_{\Of\otimes_{\mathbb{Z}}R}((\Of\otimes_{\mathbb{Z}}R)^{2})\mid
gg^{\ast}\in R^{\times}\}$$
for every $\mathbb{Z}$-algebra $R$. 
Then $\Ms^{DP}$ is a 
$G$-torsor ([14] Lemma 1.3), hence the natural forgetful functor
$$\pi:\Ms^{DP}\rightarrow\M^{DP}$$ 
is surjective and smooth. 
Consider the Grassmannian
$\Grass^{g}_{\mathbb{Z}_{(p)}}(\Of^{2}\otimes_{\mathbb{Z}}\mathbb{Z}_{(p)})$
and the closed subscheme $\N$
given by
\begin{align}
\N(R)=\{\F\subset \Of^{2}\otimes_{\mathbb{Z}}R\mid\F\text{ is an
  }\Of\otimes_{\mathbb{Z}}R\text{-submodule, locally on }\Spec R&\notag\\
\text{ a direct summand as an }R\text{-module and rk}_{R}(\F)=g&\}\notag
\end{align}
for every $\mathbb{Z}_{(p)}$-algebra $R$. Let
$\langle\cdot,\cdot\rangle$ be the 
canonical alternating $\Of$-bilinear form on $\Of^{2}$ and denote by
$\N^{DP}$
the closed subscheme of $\N$ given by
$$\N^{DP}(R)=\{\F\in \N(R)\mid \F=\F^{\perp}\},$$
where $\F^{\perp}$ is the orthogonal 
complement of $\F$ with respect to $\langle\cdot,\cdot\rangle$.
Since the module $\Phi^{\ast}(\Lie(A\d)\d)$ is totally isotropic with
respect to 
$\langle\cdot,\cdot\rangle$, we obtain a morphism
\begin{align}
f:\Ms^{DP}&\rightarrow\N^{DP}\notag\\
(A,\iota,\varphi,\eta,\Phi)&\mapsto\Phi^{\ast}(\Lie(A\d)\d).\notag
\end{align}
Using the theorem of Grothendieck-Messing ([9] 4.1.10) one can show that $f$ is
smooth.
We obtain the following diagram
$$\xymatrix{
&\ \ \Ms^{DP}\ar@{->}[dl]_\pi\ar@{->}[dr]^f&\\
\M^{DP}& &\N^{DP}.\\
}$$
Locally for the \'etale topology $\N^{DP}$ coincides with $\M^{DP}$ ([5]
Theorem 3.3). Therefore, to determine the singularities of $\M^{DP}$
it is sufficient to examine the local model. The local model is much
easier to handle 
because it can be defined in terms of linear algebra.\\

\noindent
\textbf{5.2.}
Consider the 
local model $\N^{K}$ defined by the determinant condition 
$$\N^{K}(R)=\{\F\in
\N(R)\mid\charpol(x,\F)=\prod_{\varphi:F\rightarrow\mathbb{\overline{Q}}}
(T-\varphi(x)) 
\text{ for all }x\in\Of\}.$$
The analogous construction as above shows that the local model of the
moduli problem 
$\mathcal{M}^{K}$ is equal to $\N^{DP}\cap\N^{K}$ because abelian
schemes in $\mathcal{M}^{K}$ are equipped with an $\Of$-linear
$p$-principal polarization 
and satisfy the determinant condition.

\addtocounter{prop}{2}
\begin{prop}
The generic fibres of $\N^{DP}$ and $\N^{K}$ are equal to the restriction
of scalars $\Res_{F/\mathbb{Q}}(\mathbb{P}^{1})$. In particular, they are
smooth over $\mathbb{Q}$.
\end{prop}

\noindent
\textsc{Proof:} 
Obviously $\Res_{F/\mathbb{Q}}(\mathbb{P}^{1})$ is contained in
$\N^{DP}$ and $\N^{K}$. 
On the other hand, let $\F$ be an element of $\N^{DP}(R)$ or $\N^{K}(R)$ for
a $\mathbb{Q}$-algebra $R$. We have to show that $\F$ is locally on $R$ a
free $F\otimes_{\mathbb{Q}} R$-module, hence we can assume that $R$ is
a field $k$ of 
characteristic zero. 

Let $F=\mathbb{Q}(\alpha)$ and $f$ be the
minimal polynomial of $\alpha$ over $\mathbb{Q}$. Consider the
decomposition $f=\prod g_{i}$ of $f$ into irreducible factors in $k[X]$.  
We obtain a decomposition $F\otimes_{\mathbb{Q}}k=\prod k_{i}$, where we
denote by $k_{i}$ the field $k[X]/(g_{i})$. This decomposition induces a
decomposition 
$$\F=\bigoplus_{i}\F_{i},$$ 
where $\F_{i}$ is a $k_{i}$-vector
space of dimension $n_{i}\leq 2$. We obtain that $\charpol(\alpha,\F)$ is
equal to $\prod g_{i}^{n_{i}}$. Therefore, $\F$ satisfies the determinant
condition if and only if all $n_{i}$ are equal to 1, hence if and only
if $\F$ is a free $F\otimes_{\mathbb{Q}}k$-module. On the other hand, 
$\F$ is totally isotropic with respect
to the canonical alternating $F\otimes_{\mathbb{Q}}k$-bilinear form on
$(F\otimes_{\mathbb{Q}}k)^{2}$ if and only if every $\F_{i}$ is
totally 
isotropic with respect to the canonical alternating form on
$k_{i}^{2}$ which means that all
$n_{i}$ are equal to 1. \hfill{$\Box$}

\begin{prop}
Let $p$ be unramified in $\Of$. Then
the local models $\N^{DP}$ and $\N^{K}$ are equal to
$\Res_{(\Of\otimes_{\mathbb{Z}}\mathbb{Z}_{(p)})/\mathbb{Z}_{(p)}}
(\mathbb{P}^{1})$ 
as subschemes of the 
Grassmannian.
\end{prop}

\noindent
\textsc{Proof:}
Using the
equality of the generic fibre, it is sufficient to prove the claim in the
case of a local ring $R$ with residue field $k$ of characteristic $p$.
Let $p\Of=\mathfrak{p}_{1}...\mathfrak{p}_{r}$ be the decomposition into
prime ideals of $p$ in $\Of$. This induces a decomposition
$\Of\otimes_{\mathbb{Z}}\mathbb{F}_{p}=\prod_{i=1}^{r}\Of/\mathfrak{p}_{i},$
where $\Of/\mathfrak{p}_{i}$ is a separable extension of
$\mathbb{F}_{p}$. The proposition is proved by the same
arguments as in the case of the generic fibre. \hfill{$\Box$}\\
 
\noindent
\textbf{Remark 5.5.} 
The equality
$\N^{DP}=\Res_{(\Of\otimes_{\mathbb{Z}}\mathbb{Z}_{(p)})/\mathbb{Z}_{(p)}}
(\mathbb{P}^{1})$
was already proved by Deligne and Pappas in [5] (Proposition 2.7,
Corollary 2.9).   

Let $\N^R$ be the local model of the moduli space $\M^R$ (Remark 3.8), i.e.,
$$\N^R(R)=\{\F\in\N(R)\mid \F\text{ a free }\Of\t R\text{-module of
  rank }1\}.$$ 
Obviously, $\N^R$ is equal to
  $\Res_{(\Of\otimes_{\mathbb{Z}}\mathbb{Z}_{(p)})/\mathbb{Z}_{(p)}}
  (\mathbb{P}^{1})$.  
  Thus $\N^R$ coincides with $\N^{DP}$ in the unramified case
  (Proposition 5.4) 
 which is the reason why the moduli spaces $\M^R$ and $\M^{DP}$
  coincide in this case.

\addtocounter{prop}{1}
\begin{thm} 
The local model $\N^{DP}$ is
equal to the local model $\N^{K}$ as subschemes of the Grassmannian.
\end{thm}

\noindent
\textsc{Proof:}
If $p$ is not equal to 2,
it is easy to see that $\N^{DP}$ is contained in $\N^{K}$. Indeed, let
$\F$ be an 
element of $\N^{DP}(R)$ and denote by $M$
the module $(\Of\otimes_{\mathbb{Z}}R)^{2}$. 
Since the determinant condition can be proved
locally on $R$, we may assume that $\F$ is a direct summand of $M$.
Let $d$ be an element of $D^{-1}$ such that $|D^{-1}/(d\Of)|$ is prime to
$p$. 
The module $\F$ is totally isotropic with respect to
$\langle\cdot,\cdot\rangle$, 
thus it
is totally isotropic with respect to the perfect alternating $R$-linear
form
$\langle\cdot,\cdot\rangle\p=tr_{F/\mathbb{Q}}\circ(d\langle\cdot,
\cdot\rangle)$  
(1.5). We obtain a diagram 
$$\xymatrix{
0\ar@{->}[r]&\F\ar@{->}[r]\ar@{->}[d]_h& M\ar@{->}[r]\ar@{->}[d]_\sim&
M/\F\ar@{->}[r]\ar@{->}[d]& 0\\ 
0\ar@{->}[r]&(M/\F)\d\ar@{->}[r]& M\d\ar@{->}[r]& \F\d\ar@{->}[r]& 0,\\
}$$
where we take the dual of each module with respect to the form $\langle\ ,\
\rangle\p$. The morphism $h$ is well defined because $\F$ is totally
isotropic with respect to $\langle\ ,\ 
\rangle\p$. Furthermore, it is an
injective morphism of direct summands of rank $g$ of $M$ and $M\d$
respectively, hence it is 
an isomorphism. Using the above diagram we can compute the characteristic
polynomial of an element $x\in\Of$ on $M$ 
\begin{align}
\charpol(x;M)&=\charpol(x;\F)\charpol(x;M/\F)\notag\\
&=\charpol(x;\F)^{2}.\tag{5.6.1}
\end{align}
Let $P$ be $\charpol(x;\F)$ and denote by $Q$ the polynomial
$\prod_{\sigma:F\rightarrow\overline{\mathbb{Q}}}(T-\sigma(x))$ in $R[X]$.
The characteristic polynomial of $x$ on $M=(\Of\otimes_{\mathbb{Z}}R)^{2}$ is
equal to $Q^{2}$,
hence we obtain the equation $P^{2}=Q^{2}$ from (5.6.1), i.e., 
\begin{align}
(P+Q)(P-Q)=0\tag{5.6.2}
\end{align}
in $R[X]$.
Since the leading coefficient of both polynomials is equal to 1, the
leading coefficient of the polynomial $P+Q$ is equal to 2. 
We now assume that $p$ is not equal to 2. Since $R$ is
a $\mathbb{Z}_{(p)}$-algebra, the prime 2 is invertible in $R$. Hence $P-Q$
is equal to zero and $\F$ is an element of $\N^K(R)$. Thus $\N^{DP}$
is contained in $\N^K$ if $p\neq 2$. 
In the case of characteristic 2 the polynomials $P$ and $Q$ can differ by
nilpotent elements. 

To prove the inclusion $\N^K\subset\N^{DP}$ in the case $p\neq 2$ and
the equality in the case $p=2$ we have to use the flatness of $\N^K$
in the first case and the flatness of both $\N^K$ and $\N^{DP}$ in the
second case. 
We will first reduce ourselves to the totally ramified case.
To start with, 
note that the determinant condition and the
condition to be totally isotropic can both be checked after faithfully flat
base change. Hence it is sufficient to prove the equality of the schemes
$\N^{DP}\otimes_{\mathbb{Z}_{(p)}}\mathbb{Z}_{p}$ and
$\N^{K}\otimes_{\mathbb{Z}_{(p)}}\mathbb{Z}_{p}$. Let
$p\Of=\mathfrak{p}_{1}^{e_{1}}...\mathfrak{p}_{r}^{e_{r}}$ be the
decomposition of $p$ into
prime ideals of $\Of$.
We obtain a decomposition 
$$\Of\otimes_{\mathbb{Z}}\mathbb{Z}_{p}=\prod_{i=1}^{r}
\widehat{\mathcal{O}}_{F,\mathfrak{p}_{i}},$$   
where $\widehat{\mathcal{O}}_{F,\mathfrak{p}_{i}}$ is the completion
of $\Of$ in 
$\mathfrak{p}_{i}$. This decomposition induces a decomposition of the
local models 
$$\N^{DP}=\prod_{i=1}^{r}\N_{i}^{DP}\ \text{ 
and }\ 
\N^{K}=\prod_{i=1}^{r}\N_{i}^{K}.$$
Here $\N_{i}^{DP}$ is the functor
$$
\N_{i}^{DP}(R)=\{\F_{i}\subset
(\widehat{\mathcal{O}}_{F,\mathfrak{p}_{i}}
\otimes_{\mathbb{Z}_{p}}R)^{2} \mid\F_{i}=\F_{i}^{\perp_{i}}\},$$
where $\F^{\perp_{i}}$ is the
orthogonal complement of $\F$ with respect to
the canonical alternating
$\widehat{\mathcal{O}}_{F,\mathfrak{p}_{i}}$-linear form on
$\widehat{\mathcal{O}}_{F,\mathfrak{p}_{i}}^{2}$,  
and
\begin{align}
\N_{i}^{K}(R)=\{\F_{i}\subset
(\widehat{\mathcal{O}}_{F,\mathfrak{p}_{i}}\otimes_{\mathbb{Z}_{p}}R)^{2}\mid
\charpol(x;\F_{i})=\prod_{\sigma:\widehat{F}_{\mathfrak{p}_{i}}\rightarrow
  \overline{\mathbb{Q}}_{p}}(T-\sigma(x))&\notag\\
\text{ for all }x\in\widehat{\mathcal{O}}_{F,\mathfrak{p}_{i}}\}&.\notag
\end{align}
It is sufficient to prove that $\N^{K}_{i}$ is equal to $\N^{DP}_{i}$
for all $i$. For
simplicity of notation we will write $\N^{K}$ and $F$ instead of
$\N_{i}^{K}$ and
$\widehat{F}_{\mathfrak{p}_{i}}$ etc.

Let $\Of^{ur}$  be the maximal unramified extension of $\mathbb{Z}_{p}$ in
$\Of$ and let $\kappa$ be its residue field. Since both conditions can
be proved after \'etale base change, 
it is sufficient to prove the equivalence of the local models over
$\Of^{ur}$. We have a decomposition
$$\Of\otimes_{\mathbb{Z}_{p}}\Of^{ur}=\prod_{i=1}^{f}\mathcal{O}_{F,i},$$
where $\mathcal{O}_{F,i}$ is totally ramified of degree $e$ over
$\Of^{ur}$. 
We again obtain a decomposition of the local models
$$\N^{DP}=\prod_{i=1}^{f}\N^{DP}_{i}\ \text{ and }\
\N^{K}=\prod_{i=1}^{f} \N_{i}^{K},$$
where we have
\begin{align}
\N_{i}^{K}(R)=\{\F_{i}\subset
(\mathcal{O}_{F,i}\otimes_{\mathbb{Z}_{p}}R)^{2}\mid
\charpol(x;\F_{i})=\prod_{\sigma:F_{i}\rightarrow\overline{\mathbb{Q}}_{p}}
(T-\sigma(x))&\notag\\
\text{ for all }x\in\mathcal{O}_{F,{i}}\}&\notag
\end{align}
and $\N^{DP}_{i}$ analogous. We have to show that $\N^{K}_{i}$ is
equal to $\N^{DP}_{i}$, therefore, we write again $\N^{K}$ instead of
$\N_{i}^{K}$ etc. 

Now $p$ is totally
ramified in $\Of$. 
Both models $\N^{DP}$ and $\N^{K}$ have the same generic fibre (Proposition
5.3) and they are both flat over $\Spec\mathbb{Z}_{p}$ ([5] Corollary
4.4 and [13] Corollary 4.3).  
Therefore, they are both equal to the flat closure of the generic fibre
$\N_{\eta}^{DP}=\N^{K}_{\eta}$ in $\N$. 
In the case $p\neq 2$
the elementary argument from the beginning of this proof shows that
$\N^{DP}$ is contained in $\N^{K}$. Thus when $p\neq 2$ we only need
the flatness of 
$\N^{K}$ to prove the theorem.
\hfill{$\Box$}\\

\noindent
\textbf{5.7.} We now want to give a direct proof of the equality
$\N^{K}=\N^{DP}$ in the 
totally ramified case of tame ramification. As we have seen in the proof of
Theorem 5.6, it is sufficient to prove the equality of
$\N^{DP}\otimes_{\mathbb{Z}_{(p)}}\mathbb{Z}_{p}$ and 
$\N^{K}\otimes_{\mathbb{Z}_{(p)}}\mathbb{Z}_{p}$.
For simplicity we will
change our notation and write $\N$ instead of
$\N\otimes_{\mathbb{Z}_{(p)}}\mathbb{Z}_{p}$. Let $\mathfrak{p}$ be the
prime ideal of $\Of$ over $p$. We write $F$ instead of
$\hat{F}_{\mathfrak{p}}$. Denote by $\Of$ the ring of integers of
$F$ and by $\pi$ the
uniformizing element, where
$e=e(\mathfrak{p})=g$ is the degree of the extension
$F/\mathbb{Q}_{p}$. Since the 
extension is tamely ramified, 
we can choose a uniformizing element $\pi$ such that $\pi^{e}$ is equal to
$up$, where $u$ is an element of $\mathbb{Z}_{p}$. After \'etale base
change to $\mathbb{Z}_{p}[\sqrt[e]{u}]$ we may assume
that the Eisenstein polynomial of $\pi$ is equal to
$X^{e}-p$.\\

\noindent
\textbf{5.8.}
\textsc{Proof of Theorem 5.6 in the tamely ramified case:}\\
It is sufficient to prove the equivalence of the two conditions in the
case of a local ring 
$R$. Since the claim is true for the generic fibre, we may assume that the
residue field $k$ is of characteristic $p$. 

Let $\F$ be an element of
$\N(R)$ and let $\overline{\F}$ be the reduction of $\F$
modulo
the maximal ideal of
$R$. It is an $\Of\otimes_{\mathbb{Z}}k=k[\pi]/\pi^{e}$-submodule of
$(k[\pi]/\pi^{e})^{2}$. Using the 
elementary divisor theorem for $k[[\pi]]$ we obtain a basis
$\overline{f}_{1},\overline{f}_{2}$ of $(k[\pi]/\pi^{e})^{2}$ such that
$\overline{\F}$ is generated by $\pi^{i}\overline{f}_{1}$ and
$\pi^{j}\overline{f}_{2}$ with $i\leq j$ and $i+j=e$. 
Let $f_{1}, f_{2}$ be an $\Of\otimes_{\mathbb{Z}}R$-basis of
$(\Of\otimes_{\mathbb{Z}}R)^{2}$ which lifts
$\overline{f}_{1},\overline{f}_{2}$. Consider the following submodules of
$(\Of\otimes_{\mathbb{Z}}R)^{2}$
\begin{align}\tag{5.8.1}
\begin{aligned}
\F_{0}&=\langle\pi^{e-1}f_{1},...,\pi^{i}f_{1},\pi^{e-1}f_{2},...,\pi^{j}f_{2}
\rangle_{R}\\ 
\F_{1}&=\langle\pi^{i-1}f_{1},...,f_{1},\pi^{j-1}f_{2},...,f_{2}\rangle_{R}.
\end{aligned}
\end{align}
We obtain an $R$-linear decomposition
$(\Of\otimes_{\mathbb{Z}}R)^{2}=\F_{0}\oplus\F_{1}$.  

Define the open subset $U$ of $\N$ as consisting of all modules which
are graphs of R-linear morphisms $h:\F_{0}\rightarrow\F_{1}$. 
Then $U$ is an open neighbourhood of $\F$. We will show that
$\N^{DP}\cap U$ is equal to $\N^{K}\cap U$. Since $\F$ was arbitrary, this
will prove the theorem.
As two perfect
alternating $\Of\otimes_{\mathbb{Z}}R$-bilinear forms only differ by a unit
of $\Of\otimes_{\mathbb{Z}}R$, we may use the alternating form
corresponding to the basis $f_{1}$ and $f_{2}$ instead of the alternating
form $\langle\cdot,\cdot\rangle$. We will again denote this form by
$\langle\cdot,\cdot\rangle$. \\ 

\noindent
\textbf{5.9.}
Let $\beta$ be the perfect symmetric $R$-bilinear form on
$(\Of\otimes_{\mathbb{Z}}R)^{2}$ defined by
$\beta(\pi^{i}f_{1},\pi^{j}f_{2})$ $=1$ for $i+j=e$ and
$\beta(\pi^{i}f_{1},\pi^{j}f_{2})=0$ otherwise.
Consider the morphism
\begin{align}
\phi:U&\rightarrow U\notag\\
\G&\mapsto\G^{\perp_{\beta}},\notag
\end{align}
where $\G^{\perp_{\beta}}$ is the orthogonal complement of $\G$ with
respect to $\beta$. As $\pi$ is self adjoint with respect to $\beta$, the
module $\G^{\perp_{\beta}}$ is again $\pi$-invariant. Furthermore, it is the
graph of an $R$-linear morphism $\F_{0}\rightarrow\F_{1}$, hence
$\phi$ is well defined. 

\addtocounter{prop}{3}
\begin{prop}
Let $R\p$ be an $R$-algebra and let
$\G$ be an element of $U(R\p)$.
Then $\phi(\G)=\G^{\perp_{\beta}}$ satisfies the determinant condition if
and only if $\G$ is totally isotropic with respect to
$\langle\cdot,\cdot\rangle$. 
Hence
(since $\phi^{2}$ is the 
identity) the morphism $\phi$ induces an isomorphism from $\N^{K}\cap U$ to 
$\N^{DP}\cap U$.
\end{prop}

\noindent
\textsc{Proof:}
If the integer $i$ is equal to zero, the module $\G$ is the graph of an
$R\p$-linear morphism 
$$h:\F_{0,R\p}=(\Of\otimes_{\mathbb{Z}}R\p)f_{1}\rightarrow\F_{1,R\p}=
(\Of\otimes_{\mathbb{Z}}R\p)f_{2}.$$ 
Then $\G$ is a free $\Of\otimes R\p$-submodule of rank 1 generated by
$f_{1}+h(f_{1})$, hence 
$\G$ is totally isotropic and the characteristic polynomial of the 
multiplication by $\pi$ on $\F$ is $X^{e}-p$.
Thus both conditions are trivial and the proposition is proved in this case.

Now suppose that $i$ is greater than zero.
Let $\G$ be the
graph of the $R\p$-linear morphism $h:\F_{0,R\p}\rightarrow\F_{1,R\p}$. Let
$C\in \Mat_{e}(R\p)$ be the corresponding matrix with respect to the basis of
$\F_{0}$ and $\F_{1}$ defined in (5.8.1). We will always work with the
$R\p$-basis of $(\Of\otimes_{\mathbb{Z}}R\p)^{2}$ induced by the
decomposition $\F_{0,R\p}\oplus\F_{1,R\p}$. Thus $\G$ is spanned by
the columns of the matrix  
$\bigl(\begin{smallmatrix}
1\\
C\\
\end{smallmatrix}\bigr)$ of size $2e\times e$. 
The matrix of the bilinear form $\beta$ in terms of our basis is equal to
$\bigl(\begin{smallmatrix}
  & T\\
T &\\
\end{smallmatrix}\bigr),$ 
where $T$ is equal to

$$T=\begin{pmatrix}
  &      &1\\
  &\a & \\
1 &      & \\
 \end{pmatrix}\in \Mat_{e}(R\p).$$ 
We obtain that the orthogonal complement $\G^{\perp_{\beta}}$ is
spanned by the columns of the matrix $\bigl(\begin{smallmatrix} 
1\\
 -T~^{t}CT\\
\end{smallmatrix}\bigr)$. Hence $\G^{\perp_{\beta}}$ is the graph 
of the $R\p$-linear morphism $h\p:\F_{0,R\p}\rightarrow\F_{1,R\p}$ 
corresponding to the matrix $C\p=-T~^{t}CT\in \Mat_{e}(R\p)$. 

Let $P$ be the matrix of the multiplication with $\pi$ on
$(\Of\otimes_{\mathbb{Z}}R\p)^{2}$  
with respect to the above basis.
We get
$$P=\begin{pmatrix}
K_{j}  &       &B_{ji}&      \\
       &K_{i}  &      &B_{ij}\\
pB_{ij}&       &K_{i} &      \\ 
       &pB_{ij}&      &K_{j}, \\
\end{pmatrix}$$
where $K_{i}$ and $K_{j}$ are nilpotent Jordan blocks of size $i$ and $j$
respectively and $B_{ij}$, resp. $B_{ji}$, is the $i\times j$-matrix, resp.
$j\times i$-matrix, with the entry 1 in its lower left corner and zero
elsewhere. 
Since $\G^{\perp_{\beta}}$ is $\pi$-stable, there exists an 
$e\times e$-matrix $A\p$ such that 
$$\begin{pmatrix}
                          1\\
                          C\p\\
                       \end{pmatrix}A\p=P\begin{pmatrix}
                                            1\\
                                            C\p\\
                                         \end{pmatrix}.$$
Then $A\p$ is the matrix of the $\pi$-multiplication on $\G^{\perp_{\beta}}$
in terms of the basis induced by the linear map $h\p$. We get
\begin{align}
A\p=\begin{pmatrix}
K_{j}&     \\
     &K_{i}\\
\end{pmatrix}+
\begin{pmatrix}
B_{ji}&\\
      &B_{ij}\\
\end{pmatrix}C\p.\tag{5.10.1}
\end{align}  
Thus $\G^{\perp_{\beta}}$ satisfies the determinant condition if and only
if the characteristic polynomial of $A\p$ is equal to $X^{e}-p$.

Our goal is to transform the determinant condition into the condition to be
totally isotropic for $\G$. The module $\G$ is generated by
$$e_{1}:=\pi^{i}f_{1}+h(\pi^{i}f_{1})$$
and
$$e_{2}:= \pi^{j}f_{2}+h(\pi^{j}f_{2})$$
as $\Of\otimes_{\mathbb{Z}}R\p$-module. Indeed, $e_{1},\pi
e_{1},...,\pi^{e-i-1}e_{1}, e_{2}, \pi e_{2},..., \pi^{e-j-1}e_{2}$
are contained in $\G$ since $\G$ is  
$\pi$-stable. 
Their projections on the direct summand $\F_{0,R\p}$ of
$(\Of\otimes_{\mathbb{Z}}R\p)^{2}$ generate  
$\F_{0,R\p}$, therefore, they form an $R\p$-basis of $\G$. 
Hence $\G$ is totally isotropic if and only if 
$\langle e_{1},e_{2}\rangle$ is equal to zero. 

Consider the matrix $C$ of the linear map $h$. 
Since $e_{1}$ is equal to $\pi^{i}f_{1}+h(\pi^{i}f_{1})$, the entries of  
line $j$ of $C$ are the coefficients of $h(\pi^{i}f_{1})$ in terms of the basis
$\pi^{i-1}f_{1},..., f_{1},$ $\pi^{j-1}f_{2}, ..., f_{2}$ of
$\F_{1,R\p}$. On the other  
hand, the entries of line $e$ of $C$ are the coefficients of 
$h(\pi^{j}f_{2})$. Therefore, we have
$$e_{1}=\pi^{i}f_{1}+h(\pi^{i}f_{1})=\pi^{i}f_{1}+(\sum_{l=1}^{i}c_{lj}
\pi^{i-l})f_{1}+(\sum_{k=1}^{j}c_{k+i,j}\pi^{j-k})f_{2}$$ 
and
$$e_{2}=\pi^{j}f_{2}+h(\pi^{j}f_{2})=\pi^{j}f_{2}+
(\sum_{l=1}^{i}c_{le}\pi^{i-l})f_{1}+(\sum_{k=1}^{j}c_{k+i,e}
\pi^{j-k})f_{2}.$$  
Thus we obtain
$$\langle
e_{1},e_{2}\rangle=[\pi^{i}+(\sum_{l=1}^{i}c_{lj}\pi^{i-l})][\pi^{j}+
(\sum_{k=1}^{j}c_{k+i,e}\pi^{j-k})]-(\sum_{k=1}^{j}c_{k+i,j}\pi^{j-k})
(\sum_{l=1}^{i}c_{le}\pi^{i-l}).$$  
Let $\varphi$ be the natural map from $R\p[X]$ to 
$R\p[X]/(X^{e}-p)=R\p[\pi]$ and 
$$g:=[X^{i}+(\sum_{l=1}^{i}c_{lj}X^{i-l})][X^{j}+(\sum_{k=1}^{j}c_{k+i,e}
X^{j-k})]-(\sum_{k=1}^{j}c_{k+i,j}X^{j-k})(\sum_{l=1}^{i}c_{le}X^{i-l}).$$ 
It is clear that $g$ is a normalized polynomial of degree $e$ and
$\varphi(g)=\langle e_{1},e_{2}\rangle$. Since the kernel 
of $\varphi$ is generated by the polynomial $X^{e}-p$, the element 
$\langle e_{1},e_{2}\rangle$ of $R\p[\pi]$ is equal to zero if and
only if $g$ is 
equal to $X^{e}-p$.
On the other hand, $g$ is the 
characteristic polynomial of the matrix
$$M:=\begin{pmatrix}
 0   & 1     &      &     &     &      &      &\\
     &\ddots &\ddots&     &     &      &      &\\ 
     &       &\ddots& 1   &     &      &      &\\
-c_{ee}&\ldots &\ldots&-c_{i+1,e}&-c_{i,e}&\ldots&\ldots&-c_{1,e}\\
     &       &      &     & 0   & 1    &      &\\
     &       &      &     &     &\ddots&\ddots&\\
     &       &      &     &     &      &\ddots& 1\\
-c_{e,j}&\ldots &\ldots&-c_{i+1,j}&-c_{i,j}&\ldots&\ldots&-c_{1,j}\\
\end{pmatrix}\text{,}$$
as one can verify directly.
But $M$ is equal to the matrix 
$A\p$ in (5.10.1)
of the multiplication of $\pi$ on $\G^{\perp_{\beta}}$. Thus $\G$ is
totally isotropic if and only if $\G^{\perp_{\beta}}$ satisfies the
determinant condition.\hfill{$\Box$}\\

\noindent
\textbf{5.11.}
\textsc{Continuation of the proof of Theorem 5.6:}\\
Since both conditions are trivial in the case of $i=0$, we may assume that
$i$ is greater than zero. Let $D^{-1}=(d)$ be the inverse ideal of the
different of 
the extension $F/\Qp$ and denote by $\langle\cdot,\cdot\rangle\p$ the
alternating 
$R\p$-bilinear form $\tr_{F/\Qp}\circ d\langle\cdot,\cdot\rangle$. Then $\G$
is totally 
isotropic with respect to $\langle\cdot,\cdot\rangle$ if and only if it is
totally isotropic 
with respect to $\langle\cdot,\cdot\rangle\p$ (1.5). Since $\pi$ is the
primitive element of the 
extension $F/\Qp$, the generator $d$ of $D^{-1}$ is equal to
$(e\pi^{e-1})^{-1}$. Thus the corresponding matrix of the form
$\langle\cdot,\cdot\rangle\p$ 
in terms of the usual basis is $\bigl(\begin{smallmatrix}
  & J\\
-~^{t}J &  \\
\end{smallmatrix}\bigr)$, where $J$ is the matrix
$$\begin{pmatrix}
  &      &   &  &      &1\\ 
  &      &   &  &\a & \\
  &      &   &1 &      & \\
  &      & -1&  &      & \\
  &\a &   &  &      & \\ 
-1&      &   &  &      & \\
\end{pmatrix}$$ 
with $i$-times the entry $(-1)$ and $j$-times $1$. 

Thus for $\G$ to be totally isotropic with respect to $\langle\ ,\
\rangle\p$ it  
is necessary and sufficient that 
$$\begin{pmatrix}
1&~^{t}C\\
\end{pmatrix}\begin{pmatrix}
  & J\\
-~^{t}J &  \\
\end{pmatrix}\begin{pmatrix}
1\\
C\\
\end{pmatrix}=0$$
which means that $~^{t}C$ is equal to $JCJ$.
As $\G^{\perp_{\beta}}$ is the graph corresponding to the matrix
$C\p=-T~^{t}CT$, this shows that $\G$ is totally
isotropic if and only if $\G^{\perp_{\beta}}$ is
totally isotropic. This is equivalent to $\G$ satisfying the determinant
condition (Proposition 5.10) and the theorem is proved.\hfill{$\Box$}

}
\end{document}